\RequirePackage{fix-cm}
\documentclass[smallextended,envcountsect,envcountsame]{svjour3}

\smartqed
\usepackage{amsmath,amssymb}
\usepackage{bm}
\usepackage{graphicx}
\usepackage{url}
\usepackage{booktabs}
\usepackage{float}
\usepackage[noend]{algpseudocode}
\usepackage{algorithmicx,algorithm}
\usepackage{setspace}
\usepackage{hyperref}
\usepackage{multirow}
\numberwithin{equation}{section}
\spnewtheorem{remark}{Remark}[section]{\bfseries}{\rmfamily}
\spnewtheorem{assumption}{Assumption}[section]{\bfseries}{\itshape}
  
\journalname{Journal of Scientific Computing}

\begin{document}
\title{Hybrid Newton--BFGS High-index Saddle Dynamics for Accelerating Saddle-point Searches}
\titlerunning{Hybrid Newton--BFGS High-index Saddle Dynamics}
\author{Jianyuan~Yin \and Lei~Zhang \and Shixue~Zhang}
\authorrunning{Jianyuan Yin et al.}
\institute{
Jianyuan Yin \at
School of Mathematical Sciences, 
Laboratory of Mathematics and Complex Systems, Ministry of Education, 
Beijing Normal University, Beijing 100875, China.
\email{\mbox{jyyin@bnu.edu.cn}}
\and
Lei Zhang (corresponding author) \at
School of Mathematical Sciences, Beijing International Center for Mathematical Research, Center for Quantitative Biology, Center for Machine Learning Research, Institute for Artificial Intelligence, Peking University, Beijing 100871, China.
\email{\mbox{zhangl@math.pku.edu.cn}}
\and
Shixue Zhang \at
School of Mathematical Sciences, Peking University, Beijing 100871, China.
\email{\mbox{zhangshixue@stu.pku.edu.cn}}
}
\date{}
\maketitle

\begin{abstract}
We propose a hybrid Newton--BFGS high-index saddle dynamics (HiSD) method for accelerating saddle-point searches. The method constructs a Newton correction on the unstable subspace from the available negative eigenpairs, while BFGS updates based on projected secant information approximate the inverse Hessian on the stable subspace. By combining Newton and quasi-Newton corrections, the method incorporates curvature information in both subspaces while preserving the ascent--descent structure of HiSD. Accounting for the evolving stable subspaces, we establish local $Q$-superlinear convergence of the full-memory version near nondegenerate saddle points under suitable assumptions. We further develop a limited-memory version for large-scale problems and provide fixed-step and line search strategies. 
Numerical experiments show local convergence behavior consistent with the theory and demonstrate acceleration on ill-conditioned problems and degenerate saddle manifolds in neural-network loss landscapes. The method also supports efficient upward and downward saddle searches for a liquid crystal model, further illustrating its broad applicability to saddle-point computation.

\end{abstract}

\section{Introduction}

Energy landscapes play a central role in a wide range of scientific problems, including phase transitions~\cite{Zhang2007,prl2010,Amit2014,zhang2016recent}, chemical reactions~\cite{E2010Transitionpath,heyden2005efficient,jacobson2017automated}, and neural-network training~\cite{NIPS2014,daneshmand2018,zhang2021embedding}.
Local minima represent stable or metastable states, while saddle points help characterize transitions between them~\cite{wales2003energy,E2010Transitionpath}.
In particular, index-1 saddle points correspond to transition states along minimum energy paths and determine the associated energy barriers~\cite{E2010Transitionpath,weinan2007simplified,henkelman2000improved}.
Efficiently locating saddle points is therefore essential for studying these transitions and the structure of complex energy landscapes, but remains challenging because of their instability and the high dimensionality of many applications.

Numerical methods for saddle point computation broadly fall into two categories: path-finding methods and surface-walking methods.
Path-finding methods, such as the string method \cite{weinan2007simplified,han2027} and the nudged elastic band method \cite{henkelman2000improved}, evolve paths between prescribed states to locate transition states.
Surface-walking methods search directly from an initial guess using local derivative information; examples include dimer methods \cite{Henkelman1999}, gentlest ascent dynamics \cite{gad,liu2023gad}, the iterative minimization formulation~\cite{gao2015iterative}, the minimax method  \cite{li2001minimax,liu2024nonmonotone}, the activation-relaxation method \cite{Cances2009}, gradient extremals \cite{hoffman1986gradient} and Newton trajectories \cite{hirsch2004reaction}.
Within this latter category, high-index saddle dynamics (HiSD) provides a unified framework for computing saddle points of a prescribed index~\cite{yin2019high}.
It tracks the unstable subspace of the Hessian and combines ascent in this subspace with descent in its orthogonal complement.
Together with upward and downward searches, HiSD enables the construction of solution landscapes consisting of critical points and their connections~\cite{YinPRL,yin2021searching}.
This framework and its extensions have been applied to liquid-crystal defect configurations~\cite{han2021solution}, crystal--quasicrystal transitions~\cite{yin2021transition,cui2025CPC}, excited states of rotating Bose--Einstein condensates~\cite{yin2024revealing}, and pattern formation in reaction-diffusion systems~\cite{WU2025reaction}.

Existing studies have established local convergence results and discretization error estimates for HiSD~\cite{luo2022convergence,zhang2022sinum}, and extended convergence analysis to degenerate saddle points on critical manifolds~\cite{2026degenerate}.
The local linear convergence rate of standard HiSD depends on the Hessian condition number and the accuracy of the computed eigendirections, with ill-conditioning slowing convergence and increasing the computational cost of saddle searches~\cite{luo2022convergence,2026degenerate}.
% Such ill-conditioning may affect both the unstable and stable subspaces, leading to substantial computational costs in repeated saddle searches.
% The spectral information already computed by HiSD motivates the development of Newton and quasi-Newton corrections tailored to these subspaces, with the aim of achieving faster local convergence.

Several strategies have been introduced to accelerate HiSD, including Barzilai--Borwein (BB) step sizes \cite{barzilai1988two,yin2019high}, and momentum acceleration, which improves the linear convergence rate through an inertial term~\cite{luo2025accelerated}.
More recently, subspace curvature-scaling HiSD (SCS-HiSD) uses the computed negative eigenpairs to scale the unstable directions \cite{yin2026subspace}, yielding a linear convergence rate independent of small-magnitude negative eigenvalues.
It can also accelerate selected stable directions by computing additional positive eigenpairs, but this can substantially increase the computational cost per iteration. In addition, the spectrum on the stable subspace is generally unknown a priori and may vary during iteration, making the number of additional eigenpairs difficult to prescribe.

In numerical optimization, Newton and quasi-Newton methods offer alternative ways to exploit curvature information and accelerate convergence. Newton’s method exploits curvature information by solving a linear system involving the Hessian \cite{nocedal2006numerical}, whereas BFGS builds an inverse-Hessian approximation from secant information \cite{broyden1970new,fletcher1970new,goldfarb1970family,shanno1970conditioning}, with its limited-memory variant suited to large-scale computation \cite{liu1989limited}.
Moreover, Newton's method and full-memory BFGS achieve local quadratic and superlinear convergence, respectively~\cite{nocedal2006numerical,broyden1973local,dennis1974characterization}.
Quasi-Newton techniques have also been applied to saddle-point searches \cite{kastner2008superlinearly} and minimax problems~\cite{asl2024jsymmetric}.
Motivated by these developments, we aim to accelerate searches for saddles of a prescribed index by incorporating curvature information into HiSD, while retaining its ascent–descent structure and compatibility with upward and downward searches. 
The negative eigenpairs tracked by HiSD naturally provide curvature information on the unstable subspace, while projected secant pairs provide curvature information on the stable subspace. These complementary sources motivate a hybrid approach.

We therefore propose a hybrid Newton--BFGS HiSD method that combines Newton and quasi-Newton directions on the unstable and stable subspaces, respectively.
The method constructs a Newton correction on the unstable subspace using the negative eigenpairs already computed by HiSD.
On the stable subspace, displacements and gradient differences are projected onto the same updated subspace to form secant pairs for BFGS updates, approximating the inverse Hessian on that subspace.
The two corrections jointly define the position update while preserving the ascent--descent structure of HiSD.
To reduce storage and computational costs for large-scale problems, we further develop a limited-memory version using L-BFGS updates on the stable subspace.
We provide both constant step sizes and line search strategies, with the latter allowing separate adjustment of the step sizes on the two subspaces.

The main contributions of this work are summarized as follows.
\begin{enumerate}
    \item We propose a hybrid Newton--BFGS HiSD method that combines a Newton correction on the unstable subspace with BFGS updates based on projected secant information on the stable subspace. We also develop a limited-memory version and provide constant step sizes and line search strategies.

    \item We establish local $Q$-superlinear convergence of the hybrid method with full-memory BFGS updates under suitable assumptions. The analysis controls the effects of the evolving stable subspaces on the projected secant pairs and the inverse-Hessian approximation, and establishes a projected inverse Dennis--Moré condition for the hybrid iteration.

    \item The nonlinear polynomial experiment exhibits local superlinear convergence behavior consistent with the theory. Experiments on an ill-conditioned modified Rosenbrock function and a degenerate saddle manifold in a neural network loss landscape demonstrate advantages in convergence speed and accuracy. Upward and downward searches for the Landau--de Gennes model further demonstrate the method's efficiency in large-scale saddle-point computation through comparisons of iteration counts and computation times.
\end{enumerate}

The remainder of this paper is organized as follows. Section~\ref{sec:preliminaries} introduces the preliminaries, Section~\ref{sec:hybrid-method} presents the hybrid method, and Section~\ref{sec:local-superlinear-convergence} provides the local convergence analysis. Numerical experiments are presented in Section~\ref{sec:numerical-experiments}, followed by concluding remarks in Section~\ref{sec:conclusions}. Proofs of the technical lemmas are given in Appendix~\ref{app:technical-proofs}.
\section{Preliminaries}
\label{sec:preliminaries}

\subsection{High-index saddle dynamics}
\label{subsec:hisd}

Let $E:\mathbb{R}^d\rightarrow\mathbb{R}$ be a smooth energy function, and denote its Hessian by $H(x):=\nabla^2E(x)$.
A critical point $x^*$ is called an index-$k$ saddle point if $\nabla E(x^*)=0$ and the Hessian $H(x^*)$ has exactly $k$ negative eigenvalues. We focus on nondegenerate saddle points, for which $H(x^*)$ has no zero eigenvalues.

For an index-$k$ saddle point $x^*$, let $V^*=\operatorname{span}\{v_1^*,\ldots,v_k^*\}$ be the unstable subspace associated with the $k$ negative eigenvalues of $H(x^*)$. 
Let $P_U$ denote the orthogonal projector onto a subspace $U$, 
and set $v^*=P_{V^*}x^*$ and $w^*=P_{(V^*)^\perp}x^*$. The saddle point $x^*=v^*+w^*$ then corresponds to a local solution $(v^*,w^*)$ of the minimax problem
\begin{equation}
\label{eqn:minimax}
\min_{w\in {V^*}^{\perp}}\max_{v\in V^*} E(v+w).
\end{equation}
Motivated by this local minimax structure, the HiSD for locating an index-$k$ saddle point of $E\in\mathcal{C}^2(\mathbb{R}^d,\mathbb{R})$ is given by~\cite{yin2019high}:
\begin{equation}\label{eqn:hisd}
\left\{
\begin{aligned}
\dot{x} & = - \beta
\left( I - \sum_{i=1}^k 2v_i v_i^\top \right)
\nabla E(x),\\
\dot{v}_i &=-\zeta
\left(I-v_iv_i^\top-\sum_{j=1}^{i-1}2v_jv_j^\top\right)
\nabla^2E(x)v_i, \quad i=1,\cdots,k,
\end{aligned}
\right.
\end{equation}
where $[v_i]_{i=1}^k$ are orthonormal vectors which approximate eigenvectors corresponding to the $k$ smallest eigenvalues of the Hessian $H(x)$.
$\beta$ and $\zeta$ are positive relaxation parameters.
$v^\top w$ represents the inner product of $v$ and $w$ in the Euclidean space $\mathbb{R}^d$.

A standard discrete implementation of \eqref{eqn:hisd} consists of a position update followed by a subspace update:
\begin{align}
&\ x^{(n+1)}
=
x^{(n)}
-
\beta^{(n)}
\left(
I-2\sum_{i=1}^k
v_i^{(n)}(v_i^{(n)})^\top
\right)
\nabla E(x^{(n)}),
\label{eq:standard-hisd-original}\\
&\left[
\left(
\lambda_i^{(n+1)},
v_i^{(n+1)}
\right)
\right]_{i=1}^k
=
\mathtt{EigenSol}
\left(
\left[
v_i^{(n)}
\right]_{i=1}^k,
\nabla^2E(x^{(n+1)})
\right).
\label{eq:eigensol}
\end{align}
Here $\mathtt{EigenSol}$ denotes an eigensolver that computes the required eigenpairs either accurately or approximately, using the previous directions as initial guesses.
A typical choice is the locally optimal block preconditioned conjugate gradient (LOBPCG) method \cite{knyazev2001toward}. 

For brevity, we use the following notation for the gradient, Hessian, and tracked unstable subspace at iteration \(n\):
$$
g_n:=\nabla E(x^{(n)}),\quad H_n:=H(x^{(n)}),
\quad V_n:=\operatorname{span}\left\{v_1^{(n)},\ldots,v_k^{(n)}\right\}.
$$
Let
$P_n^-:=P_{V_n}=\sum_{i=1}^k v_i^{(n)}(v_i^{(n)})^\top,
\ P_n^+:=I-P_n^-$
be the orthogonal projectors onto the tracked unstable subspace and its orthogonal complement, respectively. 
Near the target saddle point, they approximate the spectral projectors associated with the negative and positive Hessian subspaces. 

With the current subspaces $V_n$ and $V_n^\perp$, the minimax formulation \eqref{eqn:minimax} underlying HiSD motivates the following two local optimization subproblems:
\begin{equation}\label{eq:local-subproblems}
\max_{d_n^-\in V_n} E\bigl(x^{(n)}+d_n^-\bigr),
\qquad
\min_{d_n^+\in V_n^\perp} E\bigl(x^{(n)}-d_n^+\bigr).
\end{equation}
For these two subproblems, standard HiSD uses the gradient ascent and descent directions respectively, i.e.,
$d_n^-=P_n^-g_n,\ d_n^+=P_n^+g_n.$
With a common step size $\beta^{(n)}$, the original position update~\eqref{eq:standard-hisd-original} takes the form
\begin{equation}\label{eq:standard-hisd}
x^{(n+1)}=x^{(n)} +\beta^{(n)}P_n^-g_n -\beta^{(n)}P_n^+g_n.
\end{equation}
\subsection{BFGS and L-BFGS methods}
\label{subsec:bfgs-lbfgs}
We next recall the inverse-BFGS update used to approximate the inverse Hessian in unconstrained optimization \cite{nocedal2006numerical}. For two consecutive iterates $x^{(n)}$ and $x^{(n+1)}$ in solving the problem: 
$$
\min_{x\in\mathbb{R}^d}f(x),
$$
let
$s_n:=x^{(n+1)}-x^{(n)},\ y_n:=\nabla f(x^{(n+1)})-\nabla f(x^{(n)})$.
Given an inverse-Hessian approximation $G_n$, the inverse-BFGS update is given by
\begin{equation}\label{eq:standard-inverse-bfgs}
G_{n+1}=\left(I-\rho_ns_ny_n^\top\right)G_n
\left(I-\rho_ny_ns_n^\top\right)+\rho_ns_ns_n^\top,
\quad 
\rho_n=\frac{1}{y_n^\top s_n}.
\end{equation}
If $G_n$ is symmetric positive definite and the correction pair satisfies $y_n^\top s_n>0$, then $G_{n+1}$ is also symmetric positive definite.

For large-scale problems, the limited-memory BFGS (L-BFGS) method avoids forming and storing the full matrix $G_n$ \cite{liu1989limited}. Let $m$ be the prescribed memory length and $\ell\leq m$ the number of currently stored correction pairs. These pairs, denoted by $\left[(s_j,y_j)\right]_{j=1}^{\ell}$, are ordered from the oldest to the latest. 

Given a vector $q\in\mathbb{R}^d$, an initial scaling parameter $\gamma_n>0$, and the stored correction pairs, the standard two-loop recursion computes a vector
$$
r=\mathtt{TwoLoop}
\left(q;\left[(s_j,y_j)\right]_{j=1}^{\ell},\gamma_n\right)
= G_nq,
$$
where \(G_n\) is the implicit L-BFGS matrix defined by the stored correction pairs and the initial scaling \cite{liu1989limited}.
The initial inverse-Hessian approximation in the recursion is taken as
$G_n^{(0)}=\gamma_n I$.
A standard choice is $\gamma_n=\frac{s_\ell^\top y_\ell}{y_\ell^\top y_\ell}$,
where $(s_\ell,y_\ell)$ is the most recent stored correction pair. 
If no correction pair is available, one may take $G_n^{(0)}=I$.

\section{Hybrid Newton--BFGS HiSD method}
\label{sec:hybrid-method}
Standard HiSD uses gradient directions for the local optimization problems~\eqref{eq:local-subproblems}.
We now incorporate curvature information to construct more effective directions.
The unstable subspace is already tracked in HiSD through the computation of the $k$ smallest Hessian eigenpairs. Therefore, the corresponding negative Hessian block can be directly approximated by the computed eigenpairs, leading to a Newton-type correction on the evolving unstable subspace. 
On the other hand, the stable subspace is typically high-dimensional, and explicitly forming or inverting the positive Hessian block is expensive. We therefore employ a projected BFGS approximation on the stable subspace.

This leads to a hybrid Newton--BFGS realization of high-index saddle dynamics: maximization on the unstable subspace is treated by an eigenpair-based Newton correction,
whereas minimization on the stable subspace is treated by a BFGS inverse-Hessian approximation. 
The resulting method preserves the ascent--descent structure of HiSD and its ability to search for saddle points of a prescribed index, while incorporating curvature information on both subspaces without Hessian evaluations on the stable subspace.

\subsection{Eigenpair-based Newton correction on the unstable subspace}
Using the notation $g_n$, $H_n$, $V_n$, and $P_n^\pm$ introduced in Section~\ref{subsec:hisd}, let $[(\lambda_i^{(n)},v_i^{(n)})]_{i=1}^k$ denote the approximate eigenpairs computed using the same subspace update procedure as in $\eqref{eq:eigensol}$, arranged in ascending order of the eigenvalues.
 
To derive a Newton correction for the local maximization problem in~\eqref{eq:local-subproblems}, 
we first consider the ideal case near the target saddle point, in which the computed eigenpairs are exact negative eigenpairs of $H_n$.
The Newton correction on $V_n$ is then determined by the projected Newton equation \cite{nocedal2006numerical}:
\begin{equation}\label{eq:unstable-projected-newton}
P_n^-H_nP_n^-p_n^-=-P_n^-g_n,
\quad  p_n^-\in V_n.
\end{equation}
Writing $p_n^-=\sum_{i=1}^k a_i^{(n)}v_i^{(n)}$ and using
$H_nv_i^{(n)}=\lambda_i^{(n)}v_i^{(n)}$, we obtain
\[
\lambda_i^{(n)}a_i^{(n)}
=-\bigl(v_i^{(n)}\bigr)^\top g_n,
\quad  i=1,\ldots,k.
\]
Since $\lambda_i^{(n)}<0$, the unique solution of \eqref{eq:unstable-projected-newton} is
\begin{equation}\label{eq:unstable-newton-solution}
p_n^-=
\sum_{i=1}^k
\left|\lambda_i^{(n)}\right|^{-1}
v_i^{(n)}\bigl(v_i^{(n)}\bigr)^\top g_n.
\end{equation}
We use absolute values here to ensure that the operator is positive definite on $V_n$, thereby preserving the ascent direction when the tracked eigenvalues include positive values, as may occur away from the target saddle point.

We adopt this eigenpair-based correction for the local maximization problem in~\eqref{eq:local-subproblems}, and the same correction is also employed and analyzed in SCS-HiSD~\cite{yin2026subspace}.
In practice, the eigenpairs are computed numerically, and small-magnitude eigenvalues may lead to excessively large corrections.
We therefore introduce the safeguarded curvatures $\alpha_i^{(n)}=\max\{|\lambda_i^{(n)}|, \varepsilon_{\mathrm{eig}}\}$, and define the Newton operator on $V_n$ by
\begin{equation}\label{eq:unstable-newton-operator}
N_n^-=\sum_{i=1}^k
\left(\alpha_i^{(n)}\right)^{-1}
v_i^{(n)}\bigl(v_i^{(n)}\bigr)^\top.
\end{equation}
The resulting unstable-subspace correction is
\begin{equation}\label{eq:unstable-newton-direction}
d_n^-=N_n^-P_n^-g_n.
\end{equation}
Since $N_n^-P_n^-=N_n^-$, the projector $P_n^-$ in \eqref{eq:unstable-newton-direction} may be omitted in the implementation.
It is retained here to emphasize that the correction acts on the tracked unstable subspace.

\begin{remark}
The projected Newton equation above provides an idealized motivation for the eigenvalue-scaled correction, whereas the practical direction $d_n^-$ uses the computed eigenpairs and safeguarded curvatures.
Thus, $d_n^-$ is an ascent direction and serves as the eigenpair-based Newton ascent component of the hybrid iteration.
\end{remark}

\subsection{Projected BFGS approximation on the stable subspace}
\label{subsec:projected-bfgs}

For the local minimization problem in~\eqref{eq:local-subproblems}, a Newton correction would require solving a linear system on the typically high-dimensional stable subspace.
To avoid this costly solve, we use a BFGS inverse-Hessian approximation based on projected secant information to construct the search direction.
Specifically, we maintain a symmetric positive definite matrix $G_n$ on the ambient space $\mathbb{R}^d$ and define the resulting stable-subspace correction $d_n^+$ by 
\begin{equation}\label{eq:stable-bfgs-direction}
d_n^+=P_n^+G_nP_n^+g_n.
\end{equation}
Here, $G_n$ is used only through its projected action $P_n^+G_nP_n^+$, which provides the BFGS inverse-Hessian approximation restricted to the stable subspace.

Given two consecutive iterates $x^{(n)}$ and $x^{(n+1)}$, we update the tracked subspace using the same procedure as in \eqref{eq:eigensol} and form the corresponding projector $P_{n+1}^+=I-\sum_{i=1}^kv_i^{(n+1)}(v_i^{(n+1)})^\top$.
The projected displacement and gradient difference are defined by
\begin{equation}\label{eq:projected-secant-pairs}
\sigma_n=P_{n+1}^+\left(x^{(n+1)}-x^{(n)}\right),
\quad 
\eta_n=P_{n+1}^+\left(\nabla E(x^{(n+1)})-\nabla E(x^{(n)})\right).
\end{equation}
Both $\sigma_n$ and $\eta_n$ belong to the updated stable subspace $V_{n+1}^\perp$, on which we impose the inverse secant condition
\begin{equation}\label{eq:projected-inverse-secant}
\left(P_{n+1}^+G_{n+1}P_{n+1}^+\right)\eta_n=\sigma_n.
\end{equation}
This condition constrains only the projected action of $G_{n+1}$ on $\eta_n$. A sufficient condition is the full-space secant equation
$$
G_{n+1}\eta_n=\sigma_n,
$$
which can be enforced by applying the standard inverse-BFGS update to the projected pair $(\sigma_n,\eta_n)$, as a natural choice for our algorithm. 

To preserve the positive definiteness of $G_{n}$ and ensure numerical stability, we accept the projected correction pair only when it satisfies the safeguarded curvature condition
\begin{equation}\label{eq:curvature-condition}
\sigma_n^\top\eta_n
>
\mathrm{tol}_{\mathrm{curv}}
\|\sigma_n\|^2.
\end{equation}
Otherwise, the pair $(\sigma_n,\eta_n)$ is discarded.

\paragraph{Full-memory BFGS implementation.}
When \eqref{eq:curvature-condition} holds, the full-memory inverse-BFGS update is
\begin{equation}\label{eq:inverse-bfgs-update}
G_{n+1}=
\left(I-\frac{\sigma_n\eta_n^\top}{\eta_n^\top\sigma_n}\right)
G_n\left(I-\frac{\eta_n\sigma_n^\top}{\eta_n^\top\sigma_n}\right)
+\frac{\sigma_n\sigma_n^\top}{\eta_n^\top\sigma_n}.
\end{equation}
If $G_n$ is positive definite and $\sigma_n^\top\eta_n>0$, then $G_{n+1}$ remains positive definite. If \eqref{eq:curvature-condition} is not satisfied, we set $G_{n+1}=G_n$.

\paragraph{L-BFGS implementation.}
In the L-BFGS implementation, we replace the standard correction pairs $(s_j,y_j)$ in Section~\ref{subsec:bfgs-lbfgs} with the projected pairs $(\sigma_j,\eta_j)$ defined in \eqref{eq:projected-secant-pairs}. And we retain the latest $m$ pairs satisfying \eqref{eq:curvature-condition}.

If $(\sigma_{j_\ast},\eta_{j_\ast})$ is the most recently stored pair, the initial scaling in the two-loop recursion is chosen as
$$
\widehat{G}_n=\gamma_n I,
\quad 
\gamma_n=
\frac{\sigma_{j_\ast}^\top\eta_{j_\ast}}
{\eta_{j_\ast}^\top\eta_{j_\ast}}.
$$
If no admissible correction pair is available, we take $\gamma_n=1$. A pair that fails \eqref{eq:curvature-condition} is discarded, and the current scaling parameter is retained.

For the current iterate, we apply the two-loop recursion to $q_n=P_n^+g_n$ using the stored projected pairs, and denote its output by
$$
r_n=
\mathtt{TwoLoop}
\left(
q_n;[(\sigma_j,\eta_j)],\gamma_n
\right).
$$
The stable-subspace correction is then defined by
$d_n^+=P_n^+r_n$.
Here the final projection ensures that $d_n^+\in\operatorname{Range}(P_n^+)=V_n^\perp$.

\subsection{Hybrid Newton--BFGS HiSD iteration}

Combining the unstable subspace Newton correction \eqref{eq:unstable-newton-direction} with the stable subspace BFGS correction \eqref{eq:stable-bfgs-direction}, we obtain the Hybrid Newton--BFGS HiSD iteration
\begin{equation}\label{eq:hybrid-newton-bfgs-hisd}
x^{(n+1)}=x^{(n)}+\beta_1^{(n)}d_n^--\beta_2^{(n)}d_n^+,
\end{equation}
where $\beta_1^{(n)}>0$ and $\beta_2^{(n)}>0$ are the step sizes associated with the unstable and stable components, respectively.

Equivalently, substituting the two subspace corrections into
\eqref{eq:hybrid-newton-bfgs-hisd} gives
\begin{equation}\label{eq:hybrid-newton-bfgs-hisd-expanded}
x^{(n+1)}
=
x^{(n)}
+
\beta_1^{(n)}N_n^-P_n^-g_n
-
\beta_2^{(n)}P_n^+G_nP_n^+g_n.
\end{equation}
The two corrections in \eqref{eq:hybrid-newton-bfgs-hisd-expanded} are constructed independently on the tracked negative subspace and its orthogonal complement. 
Hence, the split form above, rather than a full-space quasi-Newton equation, is the fundamental operator representation of the proposed method.

\subsection{Step size selection}

The hybrid Newton--BFGS HiSD iteration \eqref{eq:hybrid-newton-bfgs-hisd} allows different step sizes for the unstable and stable corrections. The simplest strategy is to use constant step sizes
\begin{equation}
\label{eqn:constantstepsize}
\beta_1^{(n)}\equiv\beta_1,\quad \beta_2^{(n)}\equiv\beta_2.
\end{equation}
This choice avoids a complicated step size selection procedure and can be effective when suitable values are chosen for the problem. The local convergence analysis in Section~\ref{sec:local-superlinear-convergence} focuses on the unit-step case $\beta_1^{(n)}=\beta_2^{(n)}=1$.

Alternatively, we select the two step sizes through a decoupled line search based on componentwise projected-gradient reduction to improve numerical stability, using a gradient-norm criterion similar to that in~\cite{asl2024jsymmetric}.
During this search, the directions $d_n^\pm$ constructed above and the projectors $P_n^\pm$ associated with the current computed subspaces at $x^{(n)}$ remain fixed, while the two step sizes are adjusted.
Here, ``decoupled'' means that each step size is reduced only when its corresponding projected-gradient reduction condition fails, while both conditions are evaluated at the same trial point in each test.

For trial step sizes $\beta_1,\beta_2>0$, define the joint trial point
$$
x_n^{\mathrm{trial}}
=
x^{(n)}+\beta_1d_n^- -\beta_2d_n^+.
$$
The gradients at $x^{(n)}$ and the trial point are projected using the same current projectors $P_n^\pm$:
$$
g_n^\pm
=
P_n^\pm\nabla E\!\left(x^{(n)}\right),
\quad
\widetilde g_n^\pm
=
P_n^\pm\nabla E\!\left(x_n^{\mathrm{trial}}\right).
$$
Given a sufficient-reduction parameter $c_{\mathrm{ls}}\in(0,1)$, we test the conditions
\begin{equation}\label{eq:decoupled-line search}
\|\widetilde g_n^-\|
\le (1-c_{\mathrm{ls}}\beta_1)\|g_n^-\|,\quad
\|\widetilde g_n^+\|
\le (1-c_{\mathrm{ls}}\beta_2)\|g_n^+\|.
\end{equation}

Starting from $\beta_i=\overline\beta_i>0$, $i=1,2$, we accept the trial step sizes when both conditions in \eqref{eq:decoupled-line search} are satisfied and set $\beta_i^{(n)}=\beta_i$.
Otherwise, the trial step sizes are updated componentwise using a backtracking factor $\tau_{\mathrm{ls}}\in(0,1)$:
$$
\beta_i\leftarrow
\begin{cases}
\tau_{\mathrm{ls}}\beta_i,
&\text{if the corresponding condition fails},\\
\beta_i,
&\text{otherwise},
\end{cases}
\qquad i=1,2,
$$
where $i=1,2$ correspond to the unstable and stable subspace conditions, respectively.
After either step size is reduced, the joint trial point and its gradient are recomputed, and both conditions are checked again.

Given a tolerance $\varepsilon_{\mathrm{ls}}>0$, if $\|g_n^\pm\|\le\varepsilon_{\mathrm{ls}}$, the corresponding reduction condition in \eqref{eq:decoupled-line search} is regarded as satisfied.
The search is limited to $N_{\mathrm{ls}}$ tests.
If no pair of trial step sizes is accepted within this limit, the final reduced pair is used after recomputing the corresponding trial point and gradient.
In this case, the reduction conditions in \eqref{eq:decoupled-line search} need not hold.

Combining the eigenpair-based Newton correction, the projected BFGS approximation, and the above step size strategy, we summarize the complete hybrid Newton--BFGS HiSD method in Algorithm~\ref{alg:hybrid-newton-bfgs-hisd}.
In the L-BFGS implementation, $G_n$ is represented implicitly, and $d_n^+$ is computed using the two-loop recursion followed by projection onto $V_n^\perp$, as described in Section~\ref{subsec:projected-bfgs}.

\begin{algorithm}[htbp]
\caption{Hybrid Newton--BFGS HiSD for an index-$k$ saddle point.}
\label{alg:hybrid-newton-bfgs-hisd}
\setstretch{1.2}
\begin{algorithmic}[1]
\Require
Initial point $x^{(0)}\in\mathbb{R}^d$; index $k$; maximum number of iterations $N_{\max}$.
\Require Initial eigenpairs $[(\lambda_i^{(0)},v_i^{(0)})]_{i=1}^k$; 
initial BFGS approximation $G_0=\gamma_0 I, \gamma_0>0$ for full-memory BFGS; memory length $m$ for L-BFGS.
\Require
$\varepsilon_{\mathrm{eig}}>0$,
$\mathrm{tol}_{\mathrm{curv}}>0$, $\mathrm{tol}_{\mathrm{stop}}>0$.

\State Set $g_0=\nabla E(x^{(0)})$,
$P_0^-=\sum_{i=1}^k v_i^{(0)}(v_i^{(0)})^\top$,
    $P_0^+=I-P_0^-$.
\For{$n=0,1,\ldots,N_{\max}-1$}
    \State Form $N_{n}^-$ by \eqref{eq:unstable-newton-operator}.
    \State Compute
    $d_n^-=N_n^-P_n^-g_n$ and
    $d_n^+=P_n^+G_nP_n^+g_n$.
    \State Choose $\beta_1^{(n)}$, $\beta_2^{(n)}$ by \eqref{eqn:constantstepsize} or \eqref{eq:decoupled-line search}.
    \State Update
    $x^{(n+1)}=x^{(n)}+\beta_1^{(n)}d_n^--\beta_2^{(n)}d_n^+$, $g_{n+1}=\nabla E(x^{(n+1)})$.
    \If{$\|g_{n+1}\|\leq \mathrm{tol}_{\mathrm{stop}}$}
        \State \Return $x^{(n+1)}$.
    \EndIf
    \State Compute
    $[(\lambda_i^{(n+1)},v_i^{(n+1)})]_{i=1}^k
    =\mathtt{EigenSol}\bigl(
    [v_i^{(n)}]_{i=1}^k,\nabla^2 E(x^{(n+1)})
    \bigr)$.
    \State Form $P_{n+1}^-=\sum_{i=1}^k v_i^{({n+1})}(v_i^{({n+1})})^\top$,
    $P_{n+1}^+=I-P_{n+1}^-$.
    \State Set
    $\sigma_n=P_{n+1}^+(x^{(n+1)}-x^{(n)})$,
    $\eta_n=P_{n+1}^+(g_{n+1}-g_n)$.
    \If{$\sigma_n^\top\eta_n>\mathrm{tol}_{\mathrm{curv}}\|\sigma_n\|^2$}
        \State Full memory BFGS: update $G_{n+1}$ by \eqref{eq:inverse-bfgs-update};
        \State L-BFGS: update the memory using $(\sigma_n,\eta_n)$.
    \Else
        \State Full memory BFGS: $G_{n+1}=G_n$;
        \State L-BFGS: retain the current memory.
    \EndIf
\EndFor
\State \Return $x^{(N_{\max})}$.
\end{algorithmic}
\end{algorithm}

\section{Local Superlinear Convergence of the Hybrid Newton--BFGS HiSD Method}
\label{sec:local-superlinear-convergence}
In this section, we establish local $Q$-superlinear convergence of the full-memory hybrid Newton--BFGS HiSD method with exact eigenpairs and unit step sizes near a nondegenerate index-$k$ saddle point.

For our hybrid method, the exact Newton correction eliminates the first-order error on the unstable subspace, so the superlinear convergence analysis centers on the BFGS approximation on the stable subspace. Unlike standard BFGS in a fixed space~\cite{nocedal2006numerical,broyden1973local,dennis1974characterization} and $J$-symmetric methods based on a prescribed variable partition~\cite{asl2024jsymmetric}, our method uses an evolving stable spectral subspace and corresponding projected secant pairs. 
Since the coupling between the evolving subspaces and the BFGS updates prevents a direct application of classical BFGS convergence theory, we analyze the coupled iteration directly to establish a projected inverse Dennis--Moré condition.

We first specify the local setting and introduce the additional notation needed for the analysis. We then state three technical lemmas and the main convergence theorem.

\subsection{Setting and assumptions}

We retain the notations $g_n$, $H_n$, $V_n$, $P_n^\pm$, $d_n^\pm$, $G_n$, $\sigma_n$, and $\eta_n$ introduced in the preceding sections. 
For the local analysis, we additionally write
\[
B_R(x^*):=\{x\in\mathbb{R}^d:\|x-x^*\|<R\}.
\]
Unless otherwise stated, $\|\cdot\|$ denotes the Euclidean norm for vectors and the induced spectral norm for matrices, while $\|\cdot\|_F$ denotes the Frobenius norm.

To make the local regularity and spectral separation precise, we impose the following assumption, which is used throughout this section.

\begin{assumption}[Local regularity and uniform spectral separation]
\label{ass:local-regularity}
The point $x^*$ is a nondegenerate index-$k$ saddle point of $E$.
There exist constants $R>0$, $M>0$, and $0<\mu\le L$ such that:
\begin{enumerate}
\item[(i)] For all $x,y\in B_R(x^*)$,
$\|H(x)-H(y)\|\le M\|x-y\|$.

\item[(ii)] For any $x\in B_R(x^*)$, the eigenvalues
$\{\lambda_i\}_{i=1}^d$ of $H(x)$ satisfy
$\lambda_1\le\cdots\le\lambda_k<0<\lambda_{k+1}\le\cdots\le\lambda_d$
and $|\lambda_i|\in[\mu,L]$ for $1\le i\le d$.
\end{enumerate}
\end{assumption}

Throughout this section, we consider the full-memory version of the method under the following standing setting: the eigensolver returns the exact negative eigenpairs, the eigenvalue safeguard is inactive, and unit step sizes are used in each iteration, that is,
\begin{equation}
\label{ass:123}
H_n v_i^{(n)} = \lambda_i^{(n)}v_i^{(n)},
\quad
\varepsilon_{\mathrm{eig}} \leq \mu,
\quad
\beta_1^{(n)}=\beta_2^{(n)}=1.
\end{equation}
Let $P^-(x)$ and $P^+(x)$ denote the negative and positive spectral projectors of \(H(x)\), respectively. Thus, in the theoretical iteration, we have
\[
P_n^\pm=P^\pm(x^{(n)}),
\quad 
V_n=\operatorname{Ran}P_n^-,
\quad 
V_n^\perp=\operatorname{Ran}P_n^+.
\]
At the saddle point, set
\[
H_*:=H(x^*),
\ 
P_*^\pm:=P^\pm(x^*),
\ 
V_*:=\operatorname{Ran}P_*^-=V^*,
\ 
V_*^\perp:=\operatorname{Ran}P_*^+.
\]
The inverse-BFGS matrix is initialized by $G_0=G_0^\top\succ0$. Whenever a nonzero projected correction pair is generated, the full-memory inverse-BFGS update \eqref{eq:inverse-bfgs-update} is applied. An update is called well defined when $\sigma_n^\top\eta_n>0$. Consequently, every well-defined sequence satisfies $G_n=G_n^\top\succ0,\ n\ge0$.
All the results in this section are understood under this standing setting unless otherwise stated.

Under the setting \eqref{ass:123}, we define the actual step by
\begin{equation}
\label{eq:local-split-step}
s_n
:=x^{(n+1)}-x^{(n)}
 =d_n^--d_n^+ 
 =
-\left(H_n|_{V_n}\right)^{-1}P_n^-g_n
-P_n^+G_nP_n^+g_n.
\end{equation}
With this notation, the projected displacement introduced in
\eqref{eq:projected-secant-pairs} is $\sigma_n=P_{n+1}^+s_n$.
The local analysis also uses the following assumption on the generated sequence.

\begin{assumption}[Projected-step condition]
\label{ass:projected-step}
There exists a constant $\chi\in(0,1]$ such that every nonzero step satisfies
\begin{equation}
\label{eq:Projected-step-condition}
\|\sigma_n\|
=
\|P_{n+1}^+s_n\|
\ge \chi\|s_n\|.
\end{equation}
\end{assumption}

Assumption~\ref{ass:projected-step} is a condition on the generated sequence rather than an additional regularity assumption on the energy function. 
In all subsequent results, the constant $\chi\in(0,1]$ is fixed in advance, and the local constants and
neighborhood sizes are allowed to depend on $\chi$.

\subsection{Positive-space error measures}

We regard every operator on $V_n^\perp$ as an operator on
$\mathbb{R}^d$ by extending it by zero on $V_n$. Define the positive square-root operator
\[
S(x):=\left(P^+(x)H(x)P^+(x)\right)^{1/2},
\quad
S_n:=S(x^{(n)}),
\]
which satisfies $\|S_n z\|\in [\sqrt{\mu}, \sqrt{L}]$ for any $z \in V_n^\perp$ with $\|z\|=1$. 
Whenever an inverse or a negative fractional power of such an operator appears, it is taken on $V_n^\perp$ and then extended by zero on $V_n$. In particular,
\[
S_n^{-1}S_n=S_nS_n^{-1}=P_n^+.
\]

The positive-space inverse error and its Frobenius norm are defined by
\begin{equation}
\label{eq:current-positive-error}
D_n
:=
S_n\left(P_n^+G_nP_n^+\right)S_n-P_n^+,
\quad 
\Phi_n:=\|D_n\|_F.
\end{equation}
Before applying the $n$th inverse-BFGS update, the corresponding quantities in the new positive subspace are
\begin{equation}
\label{eq:preupdate-positive-error}
\widehat D_n
:=
S_{n+1}
\left(P_{n+1}^+G_nP_{n+1}^+\right)
S_{n+1}
-P_{n+1}^+,
\quad 
\widehat\Phi_n:=\|\widehat D_n\|_F.
\end{equation}
In particular, $D_n=0$ precisely when
$\left.
\left(P_n^+G_nP_n^+\right)
\right|_{V_n^\perp}
=
\left(H_n|_{V_n^\perp}\right)^{-1}$.

The main theorem will be stated using an initial error measured on the fixed positive subspace at the saddle point:
\begin{equation}
\label{eq:initial-positive-error}
\Delta_0^+
:=
\left\|
\left.
\left(P_*^+G_0P_*^+\right)
\right|_{V_*^\perp}
-
\left(H_*|_{V_*^\perp}\right)^{-1}
\right\|_F.
\end{equation}
Thus, $\Delta_0^+$ measures the initial inverse-Hessian approximation error on the fixed positive subspace at $x^*$, whereas $\Phi_n$ and $\widehat\Phi_n$ measure the weighted approximation errors on the current and updated positive subspaces, respectively.

\subsection{Technical lemmas}
We first state the local perturbation estimates needed to compare two consecutive positive spectral subspaces and to control the projected secant pairs. The proofs of the three lemmas below are given in Appendix~\ref{app:technical-proofs}.

\begin{lemma}[Local spectral, secant, and moving-space estimates]
\label{lem:local-moving-estimates}
Under the setting \eqref{ass:123}, suppose that Assumptions~\ref{ass:local-regularity} and
\ref{ass:projected-step} hold, then there exist
$R_0\in(0,R]$ and $C_{\rm loc}>0$, with $C_{\rm loc}$ depending only on
$M$, $\mu$, $L$, and $d$, such that, whenever
$x^{(n)},x^{(n+1)}\in B_{R_0}(x^*)$,
\begin{align}
&\|P_{n+1}^\pm-P_n^\pm\|_F
+\|S_{n+1}-S_n\|_F
\le C_{\rm loc}\|s_n\|,
\label{eq:projector-square-root-estimate}\\
&\|\eta_n-H_{n+1}\sigma_n\|
\le \frac{M}{2}\|s_n\|^2,
\label{eq:projected-secant-estimate}\\
&\|\widehat D_n- D_n\|_F
\le
C_{\rm loc}\bigl(1+\|G_n\|\bigr)\|s_n\|.
\label{eq:moving-potential-comparison}
\end{align}
If $s_n\ne0$, then
\begin{equation}
\label{eq:local-curvature-lower-bound}
\sigma_n^\top\eta_n
\ge
\left(
\mu-\frac{M}{2\chi}\|s_n\|
\right)\|\sigma_n\|^2.
\end{equation}
In particular, $\sigma_n^\top\eta_n \ge
\frac{\mu}{2}\|\sigma_n\|^2>0$ whenever $\|s_n\|\le\mu\chi/M$.

Furthermore, for every $\Gamma_0>0$, if $x^{(0)}\in B_{R_0}(x^*)$ and $\|G_0\|\le\Gamma_0$, then
\begin{equation}
\label{eq:initial-potential-comparison}
\Phi_0 \le L\Delta_0^+ +C_{\rm loc}(1+\Gamma_0)\|x^{(0)}-x^*\|.
\end{equation}
\end{lemma}

The next lemma adapts the classical bounded-deterioration approach in quasi-Newton analysis~\cite{broyden1973local} to the full-memory inverse-BFGS update with evolving positive spectral subspaces.
It also controls the part of the full-space matrix $G_n$ that is not measured directly by $\Phi_n$.

For a nonzero step, define the normalized weighted secant direction by
\begin{equation}
\label{eq:weighted-secant-direction}
\omega_n
:=
\frac{S_{n+1}\sigma_n}
     {\|S_{n+1}\sigma_n\|}.
\end{equation}

\begin{lemma}[Stability and bounded deterioration of the BFGS matrices]
\label{lem:bounded-deterioration}
Under the setting \eqref{ass:123}, suppose that Assumptions~\ref{ass:local-regularity} and
\ref{ass:projected-step} hold, and fix
$\overline\Phi\in(0,1)$.
\begin{enumerate}
\item[(i)] \emph{Uniform boundedness of $G_n$.}
For every $\Lambda>0$ and $\Gamma_0>0$, there exists a constant $C_G>0$, independent of $N$, 
% if a finite sequence $\{(x^{(j)},G_j)\}_{j=0}^N$ of well-defined inverse-BFGS updates satisfies 
such that, if any finite sequence $\{(x^{(j)},G_j)\}_{j=0}^N$ generated by well-defined inverse-BFGS updates satisfies $x^{(j)}\in B_{R_0}(x^*),\ 0\le j\le N$, and
\[
G_0=G_0^\top\succ0,
\
\|G_0\|\le\Gamma_0,
\ 
\max_{0\le j\le N}\Phi_j\le\overline\Phi,
\
\sum_{j=0}^{N-1}\|s_j\|\le\Lambda,
\]
% then we have a uniform bound of $\|G_j\|$, i.e.,
then we have the uniform bound property
\begin{equation}
\label{eq:finite-prefix-G-bound}
\max_{0\le j\le N}\|G_j\|\le C_G.
\end{equation}

\item[(ii)] \emph{One-step bounded deterioration.}
For every $\overline\Gamma>0$, there exist constants
$h>0$, $c_{\rm B}>0$, and $C_{\rm B}>0$ such that, if
$x^{(n)},x^{(n+1)}\in B_{R_0}(x^*)$, and
\[
G_n=G_n^\top\succ0,
\ 
\|G_n\|\le\overline\Gamma,
\ 
\max\{\Phi_n,\widehat\Phi_n\}\le\overline\Phi,
\ 
0<\|s_n\|\le h,
\]
then the inverse-BFGS update \eqref{eq:inverse-bfgs-update} is well defined,
$G_{n+1}=G_{n+1}^\top\succ0$, and
\begin{equation}
\label{eq:bfgs-deterioration-current-space}
\Phi_{n+1}
\le
\Phi_n
-c_{\rm B}\|\widehat D_n\omega_n\|^2
+C_{\rm B}\|s_n\|.
\end{equation}

\end{enumerate}
\end{lemma}

Under uniform boundedness and step-summability assumptions, the final lemma establishes a projected inverse analogue of the classical Dennis--Mor\'e condition~\cite{dennis1974characterization} on the current positive spectral subspace.
This condition, together with the exact Newton correction on the unstable subspace, yields a vanishing relative Newton residual.
The theoretical iteration is said to terminate when $g_n=0$. A nonterminating sequence therefore satisfies $g_n\ne0$ for every $n$.

\begin{lemma}[Projected inverse Dennis--Mor\'e condition and Newton residual]
\label{lem:inverse-dennis-more}
Under the setting \eqref{ass:123}, let $\{x^{(n)}\}$ be a nonterminating sequence generated by the unit-step full-memory iteration~\eqref{eq:local-split-step}.
Suppose that all iterates lie in $B_{R_0}(x^*)$, all inverse-BFGS updates are well defined, Assumptions \ref{ass:local-regularity} and \ref{ass:projected-step} hold, and, for some $\overline\Gamma>0$ and $\overline\Phi\in(0,1)$, we assume
\begin{equation}
\sup_{n\ge0}\|G_n\|\le\overline\Gamma,
\quad 
\sup_{n\ge0}\max\{\Phi_n,\widehat\Phi_n\}
\le\overline\Phi,
\quad 
\sum_{n=0}^{\infty}\|s_n\|<\infty.
\end{equation}
Then we have
\begin{equation}
\label{eq:current-space-inverse-dennis-more}
\frac{
\left\|
\left(
P_n^+G_nP_n^+H_nP_n^+-P_n^+
\right)P_n^+s_n
\right\|
}{
\|s_n\|
}
\longrightarrow0, \quad\text{as} \ n\rightarrow \infty.
\end{equation}

Consequently, the Newton residual $\mathcal R_n:=g_n+H_ns_n$ satisfies
\begin{equation}
\label{eq:newton-residual}
P_n^-\mathcal R_n=0,
\quad
\frac{\|\mathcal R_n\|}{\|s_n\|}
\longrightarrow0, \quad\text{as} \  n\rightarrow \infty.
\end{equation}
\end{lemma}
\subsection{Main convergence theorem}

We are now ready to state the local convergence result. Its initialization conditions are expressed in terms of $\Delta_0^+$ defined in $\eqref{eq:initial-positive-error}$. 

\begin{theorem}[Local linear and Q-superlinear convergence]
\label{thm:local-superlinear-convergence}
Under the setting \eqref{ass:123}, suppose that Assumptions~\ref{ass:local-regularity} and
\ref{ass:projected-step} hold with constant $\chi\in(0,1]$, and fix $\Gamma_0>0$, $r\in(0,1)$. Then there exist sufficiently small $\varepsilon\in(0,R_0)$ and $\delta>0$ such that, if
\[
G_0=G_0^\top\succ0,
\quad 
\|G_0\|\le\Gamma_0,
\quad 
\|x^{(0)}-x^*\|\le\varepsilon,
\quad 
\Delta_0^+\le\delta,
\]
then the full-memory iteration is well defined,
remains in $B_{R_0}(x^*)$, and converges to $x^*$. Moreover, we have
\begin{equation}
\label{eq:prescribed-linear-contraction}
\|x^{(n+1)}-x^*\|
\le
r\|x^{(n)}-x^*\|
\end{equation}
for every $n$ before termination. The method either reaches $x^*$ in
finitely many iterations or satisfies
\begin{equation}
\label{eq:q-superlinear-convergence}
\lim_{n\to\infty}
\frac{\|x^{(n+1)}-x^*\|}
     {\|x^{(n)}-x^*\|}
=0.
\end{equation}
\end{theorem}

\begin{proof}
Write $e_n=x^{(n)}-x^*$ and $K_n=P_n^+G_nP_n^+$.
By the initial comparison \eqref{eq:initial-potential-comparison} in Lemma~\ref{lem:local-moving-estimates}, we have
$\Phi_0 \le L\Delta_0^+ +C_{\rm loc}(1+\Gamma_0)\|x^{(0)}-x^*\|$.
Thus, sufficiently small initial position and matrix errors ensure that $\Phi_0$ is small.
We first derive a one-step error estimate and then establish the linear contraction together with uniform bounds on $G_n$ and $\Phi_n$ by induction.

By the definition of $D_n$ in \eqref{eq:current-positive-error}, we have
\[
K_n=S_n^{-1}(P_n^++D_n)S_n^{-1},
\quad 
K_nH_nP_n^+-P_n^+=S_n^{-1}D_nS_n.
\]
Since $\|S_n^{-1}|_{V_n^{\perp}}\|\le\mu^{-1/2}$ and $\|D_n\|\le\Phi_n$, whenever $\Phi_n\le\overline\Phi<1$, we have
\begin{equation}
\label{eq:Kn-bound}
\|K_n\|\le\|S_n^{-1}|_{V_n^{\perp}}\|^2(1+\|D_n\|)\le\frac{1+\Phi_n}{\mu}\le\frac{1+\overline\Phi}{\mu}.
\end{equation}
For any $x^{(n)}\in B_{R_0}(x^*)$, the fundamental theorem of calculus and $\nabla E(x^*)=0$ give
\[
g_n-H_ne_n
=\int_0^1\bigl[H(x^*+te_n)-H_n\bigr]e_n\,dt.
\]
The Lipschitz continuity of the Hessian in Assumption~\ref{ass:local-regularity} yields
\begin{equation}\label{eq:main-proof-taylor-remainder}
\|g_n-H_ne_n\|
\le\int_0^1M(1-t)\|e_n\|^2\,dt
=\frac M2\|e_n\|^2.
\end{equation}
Since $(H_n|_{V_n})^{-1}P_n^-H_n=P_n^-$ and $K_nH_n=K_nH_nP_n^+$, the step formula \eqref{eq:local-split-step} gives
\[
\begin{aligned}
e_{n+1}
&=e_n-\left[(H_n|_{V_n})^{-1}P_n^-+K_n\right]g_n\\
&=(P_n^+-K_nH_nP_n^+)e_n
-\left[(H_n|_{V_n})^{-1}P_n^-+K_n\right](g_n-H_ne_n)\\
&=-S_n^{-1}D_nS_ne_n
-\left[(H_n|_{V_n})^{-1}P_n^-+K_n\right](g_n-H_ne_n).
\end{aligned}
\]
$(H_n|_{V_n})^{-1}P_n^-$ and $K_n$ act on $V_n$ and $V_n^\perp$ respectively, so
$\|(H_n|_{V_n})^{-1}P_n^-\| \le 1/\mu$ and \eqref{eq:Kn-bound} indicate
$\|(H_n|_{V_n})^{-1}P_n^-+K_n\|\le(1+\overline\Phi)/\mu$.
Hence, together with \eqref{eq:main-proof-taylor-remainder}, we have
\begin{equation}\label{eq:main-proof-one-step}
\|e_{n+1}\|
\le
\sqrt{\frac L\mu}\,\Phi_n \|e_n\|
+\frac{M(1+\overline\Phi)}{2\mu}\|e_n\|^2.
\end{equation}
Estimate \eqref{eq:main-proof-one-step} requires only $x^{(n)}\in B_{R_0}(x^*)$ and $\Phi_n\le\overline\Phi$; it can subsequently be used to establish that $x^{(n+1)}\in B_{R_0}(x^*)$.

Choose $\overline\Phi\in(0,1)$ such that $\sqrt{L/\mu}\,\overline\Phi\le r/2$.
Let $C_G$ be the constant given by Lemma~\ref{lem:bounded-deterioration}(i)
for $\Gamma_0$, $\overline\Phi$, and $\Lambda=R_0$, and set
$\overline\Gamma=\max\{\Gamma_0,C_G\}$.
For $\overline\Gamma$ and $\overline\Phi$, let $h,c_{\rm B},C_{\rm B}>0$ be the constants in
Lemma~\ref{lem:bounded-deterioration}(ii).

Our aim is to prove by induction that the iteration is well defined before termination and satisfies
\begin{equation}\label{eq:main-proof-induction}
\begin{gathered}
\|e_j\|\le r^j\|e_0\|,
\  G_j=G_j^\top\succ0,
\  \|G_j\|\le\overline\Gamma,\\
\Phi_j\le\Phi_0+C_{\rm B}\sum_{\ell=0}^{j-1}\|s_\ell\|
\le\frac{\overline\Phi}{2}.
\end{gathered}
\end{equation}
To establish \eqref{eq:main-proof-induction}, choose $\varepsilon\in(0,R_0)$ and $\delta>0$ sufficiently small such that
\begin{equation}\label{eq:main-proof-smallness}
\begin{gathered}
\frac{M(1+\overline\Phi)}{2\mu}\varepsilon\le\frac r2,
\quad 
(1+r)\varepsilon\le h,
\quad 
\frac{1+r}{1-r}\varepsilon\le R_0,\\
C_{\rm loc}(1+\overline\Gamma)(1+r)\varepsilon
\le\frac{\overline\Phi}{2},\\
L\delta+\left[C_{\rm loc}(1+\Gamma_0)
+C_{\rm B}\frac{1+r}{1-r}\right]\varepsilon
\le\frac{\overline\Phi}{2}.
\end{gathered}
\end{equation}
The first three conditions in \eqref{eq:main-proof-smallness} will ensure contraction, sufficiently small steps, and a bound on the total step length.
The last two conditions will control the moving-space error and the accumulated increase in $\Phi_n$.

The first condition in \eqref{eq:main-proof-smallness} implies $\varepsilon<\mu/M$.
If the iterate $x^{(n)}\in\overline{B_\varepsilon(x^*)}$, then \eqref{eq:main-proof-taylor-remainder} and the spectral bound $\|H_ne_n\|\ge\mu\|e_n\|$ give
\begin{equation}\label{eq:main-proof-gradient-lower-bound}
\|g_n\|
\ge\|H_ne_n\|-\|g_n-H_ne_n\|
\ge\left(\mu-\frac M2\varepsilon\right)\|e_n\|
\ge\frac\mu2\|e_n\|.
\end{equation}
Thus a stationary iterate in $\overline{B_\varepsilon(x^*)}$ must be $x^*$.

We now verify \eqref{eq:main-proof-induction} for $j=0$.
If $g_0=0$, \eqref{eq:main-proof-gradient-lower-bound} gives $x^{(0)}=x^*$, and the conclusion follows.
Otherwise, \eqref{eq:initial-potential-comparison} and the last inequality in \eqref{eq:main-proof-smallness} give
\[
\Phi_0\le L\delta+C_{\rm loc}(1+\Gamma_0)\varepsilon
\le\frac{\overline\Phi}{2}.
\]
Together with $G_0=G_0^\top\succ0$ and $\|G_0\|\le\Gamma_0\le\overline\Gamma$, this proves \eqref{eq:main-proof-induction} for $j=0$.

Suppose that \eqref{eq:main-proof-induction} holds for all $0\le j\le n$ and that $g_n\ne0$.
Positive definiteness of $G_n$ and the step formula \eqref{eq:local-split-step} imply $s_n\ne0$.
First compute $s_n$ and $x^{(n+1)}=x^{(n)}+s_n$ from $x^{(n)},G_n$.
By \eqref{eq:main-proof-one-step}, the induction hypothesis \eqref{eq:main-proof-induction}, and \eqref{eq:main-proof-smallness},
\[
\|e_{n+1}\|
\le\left(\sqrt{\frac L\mu}\,\overline\Phi
+\frac{M(1+\overline\Phi)}{2\mu}\varepsilon\right)\|e_n\|
\le r\|e_n\|.
\]
Thus, we have    $\|e_{n+1}\|\le\varepsilon<R_0$.
If $g_{n+1}=0$, \eqref{eq:main-proof-gradient-lower-bound} gives $x^{(n+1)}=x^*$, and the method terminates at the saddle point.
We therefore assume $g_{n+1}\ne0$ for the remainder of the induction step.

Moreover,
$\|s_n\|\le\|e_n\|+\|e_{n+1}\|
\le(1+r)r^n\|e_0\|\le h$.
Applying $\|s_j\|\le(1+r)r^j\|e_0\|$ for $0\le j\le n$, we obtain
\[
\sum_{j=0}^n\|s_j\|
\le(1+r)\|e_0\|\sum_{j=0}^nr^j
\le\frac{1+r}{1-r}\varepsilon\le R_0.
\]

We next establish the bounds on $G_{n+1}$ and $\Phi_{n+1}$ in \eqref{eq:main-proof-induction}.
Since $x^{(n)},x^{(n+1)}\in B_{R_0}(x^*)$, the moving-space estimate \eqref{eq:moving-potential-comparison}, the bounds $\|G_n\|\le\overline\Gamma$ and $\Phi_n\le\overline\Phi/2$, and \eqref{eq:main-proof-smallness} yield
\[
\widehat\Phi_n
\le\Phi_n+\|\widehat D_n-D_n\|_F
\le\frac{\overline\Phi}{2}
+C_{\rm loc}(1+\overline\Gamma)\|s_n\|
\le\overline\Phi.
\]
From $s_n\ne0$, Assumption~\ref{ass:projected-step} also gives $\|\sigma_n\|\ge\chi\|s_n\|>0$. Together with the bounds established above, this verifies all hypotheses of Lemma~\ref{lem:bounded-deterioration}(ii).
The $n$th update is well defined and $G_{n+1}=G_{n+1}^\top\succ0$. Dropping the nonpositive term in \eqref{eq:bfgs-deterioration-current-space} and using \eqref{eq:main-proof-induction}, we obtain
\[
\begin{aligned}
\Phi_{n+1}
&\le\Phi_n+C_{\rm B}\|s_n\|
\le\Phi_0+C_{\rm B}\sum_{j=0}^n\|s_j\|\\
&\le L\delta+\left[C_{\rm loc}(1+\Gamma_0)+C_{\rm B}\frac{1+r}{1-r}\right]\varepsilon\le\frac{\overline\Phi}{2}.
\end{aligned}
\]
Here the last line follows from \eqref{eq:initial-potential-comparison}, $\sum_{j=0}^n\|s_j\|\le(1+r)\varepsilon/(1-r)$, and the last inequality in \eqref{eq:main-proof-smallness}.

The first $n+1$ updates are now well defined. Moreover, $x^{(j)}\in B_{R_0}(x^*)$ and $\Phi_j\le\overline\Phi$ for $0\le j\le n+1$, while $\sum_{j=0}^{n}\|s_j\|\le R_0$.
Applying Lemma~\ref{lem:bounded-deterioration}(i) with $N=n+1$ and $\Lambda=R_0$ therefore yields $\|G_{n+1}\|\le C_G\le\overline\Gamma$, which completes the induction.

Consequently, every generated iterate remains in $B_{R_0}(x^*)$, and
\eqref{eq:prescribed-linear-contraction} holds at every step before termination.
By \eqref{eq:main-proof-gradient-lower-bound}, finite termination can occur only at $x^*$.
It remains to consider a nonterminating sequence. The contraction \eqref{eq:prescribed-linear-contraction} gives $x^{(n)}\to x^*$ and
\[
\sum_{n=0}^{\infty}\|s_n\|\le\frac{1+r}{1-r}\varepsilon<\infty.
\]
The induction also establishes that all updates are well defined,
$\sup_n\|G_n\|\le\overline\Gamma$, and $\sup_n\max\{\Phi_n,\widehat\Phi_n\}\le\overline\Phi<1$. All hypotheses of Lemma~\ref{lem:inverse-dennis-more} therefore hold, and hence $\|\mathcal R_n\|/\|s_n\|\to0$.

Finally, the identity $\mathcal R_n=g_n+H_ns_n$ gives $e_{n+1}=H_n^{-1}(\mathcal R_n+H_ne_n-g_n)$.
Using $\|H_n^{-1}\|\le1/\mu$, \eqref{eq:main-proof-taylor-remainder}, and $\|s_n\|\le(1+r)\|e_n\|$, we obtain
\[
% \begin{aligned}
\frac{\|e_{n+1}\|}{\|e_n\|}
\le\frac{\|\mathcal R_n\|}{\mu\|e_n\|}
+\frac{M}{2\mu}\|e_n\|
\le\frac{1+r}{\mu}\frac{\|\mathcal R_n\|}{\|s_n\|}
+\frac{M}{2\mu}\|e_n\|\longrightarrow0, \quad\text{as} \  n\rightarrow \infty.
% \end{aligned}
\]
This proves \eqref{eq:q-superlinear-convergence}.\qed
\end{proof}

\section{Numerical Experiments}
\label{sec:numerical-experiments}

In this section, we present four numerical experiments to examine the performance of the hybrid Newton--BFGS HiSD method.
We first consider a nonlinear polynomial with a known saddle point to examine the local convergence behavior of the full-memory and limited-memory versions.
We then compare different methods on an ill-conditioned modified Rosenbrock problem.
Next, we consider a neural-network loss function to examine the acceleration achieved near a degenerate saddle manifold.
Finally, we apply the limited-memory version to the Landau--de Gennes model and investigate its performance in upward and downward saddle searches.

\subsection{Nonlinear polynomials: local convergence behavior}

We first consider the nonlinear polynomial energy
\[
E(x)
=
\sum_{i=1}^{k}
\left(
-\frac12 a_i x_i^2+\frac{\gamma}{4}x_i^4
\right)
+
\sum_{j=1}^{d-k}
\frac12 b_j
\left(
x_{k+j}-\sum_{i=1}^{k}\alpha_{ji}x_i^2
\right)^2,
\]
where $x=(x_1,\ldots,x_d)^\top\in\mathbb R^d$ and $a_i,b_j>0$.
The origin $x^*=0$ is a nondegenerate index-$k$ saddle point. The nonlinear coupling makes the Hessian subspaces vary with the iteration point, allowing us to examine convergence in the presence of evolving subspaces.

In our experiment, we take $d=60$, $k=6$, and $\gamma=0.2$.
We set the coefficients \(a_i=1+\frac{100(i-1)}{k-1}\), \(i=1,\ldots,k\); \(b_j=1+\frac{100(j-1)}{d-k-1}\), \(j=1,\ldots,d-k\), and $\alpha_{ji}$ are sampled uniformly from $[0.5,20.5]$.
Each component of the initial point is sampled uniformly from $[-5\times10^{-3},5\times10^{-3}]$.
We compare the full-memory and limited-memory line search versions, with memory length $m=5$ for the latter.
Both methods start from the same initial point, use accurate Hessian eigenpairs at each iteration, and initialize each line search with trial step sizes $\overline{\beta}_1=\overline{\beta}_2=1$.
The iterations terminate when the gradient norm is at most $10^{-11}$.

\begin{figure}[!htbp]
    \centering
    \includegraphics[width=1\textwidth]{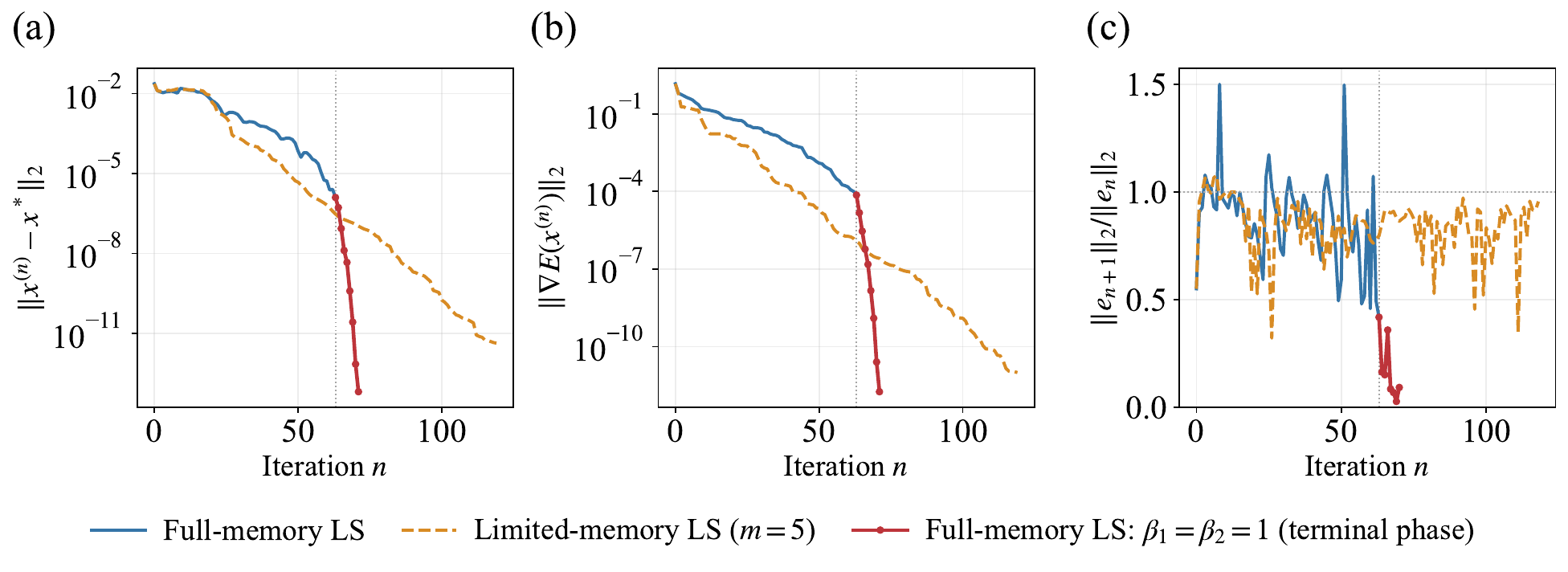}
    \caption{Plots of (a) error $\|e_n\|_2$, (b) gradient norm $\|\nabla E(x^{(n)})\|_2$, and
    (c) successive error ratio $\|e_{n+1}\|_2/\|e_n\|_2$ versus iteration number for the full-memory and limited-memory line search versions on the nonlinear polynomial example.
    The red segments indicate the terminal unit-step phase ($\beta_1^{(n)}=\beta_2^{(n)}=1$) of the full-memory version.}
    \label{fig:nonlinear-polynomial}
\end{figure}

Since the target saddle is the origin, we can directly evaluate the error norm $\|e_n\|_2=\|x^{(n)}\|_2$ and the successive error ratio $\|e_{n+1}\|_2/\|e_n\|_2$. Figure~\ref{fig:nonlinear-polynomial} shows the error norm, gradient norm, and error ratio versus iteration number.
The red segments indicate the terminal phase in which the full-memory version continuously takes unit steps, $\beta_1^{(n)}=\beta_2^{(n)}=1$.

As illustrated in Figure~\ref{fig:nonlinear-polynomial}, both methods approach the target saddle, but exhibit different terminal behaviors.
The limited-memory error curve is approximately straight on the semi-logarithmic plot, while its error ratio fluctuates within a range.
In contrast, the full-memory version exhibits pronounced acceleration after entering the consecutive unit-step phase.
Both the error and gradient norm decrease rapidly, and the error ratio becomes substantially smaller.
Although the error ratio is not monotonically decreasing, this terminal behavior is consistent with the local Q-superlinear convergence result for the full-memory version.

\subsection{Modified Rosenbrock function: an ill-conditioned problem}
We next consider the modified Rosenbrock function used in~\cite{luo2025accelerated} to evaluate the methods on an ill-conditioned saddle-point problem.
The construction starts from the $d$-dimensional Rosenbrock function
\[
R(x)
=
\sum_{i=1}^{d-1}
\left[
100(x_{i+1}-x_i^2)^2+(1-x_i)^2
\right],
\]
whose minimizer is $x^*=(1,\ldots,1)^\top$.
The modified function is given by
\[
R_m(x)
=
R(x)+\sum_{i=1}^{d}s_i\arctan^2(x_i-x_i^*).
\]
In our experiment, we set $d=1000$ and
$s_i=
\begin{cases}
-50000, & 1\le i\le5,\\
1, & 6\le i\le d,
\end{cases}$
and $x^*$ is an index-5 saddle point.
The ratio between the largest and smallest Hessian eigenvalue
magnitudes at $x^*$ is approximately $3.99\times10^4$,
indicating a substantial difference in curvature scales.

We compare standard HiSD \cite{yin2019high}, HiSD with Barzilai--Borwein step sizes (HiSD-BB) \cite{yin2019high}, accelerated HiSD (A-HiSD) \cite{luo2025accelerated}, and three versions of the proposed method: full-memory with line search, limited-memory with fixed step sizes, and limited-memory with line search.
The initial point is obtained by applying a random perturbation of Euclidean norm $0.1$ to $x^*$. All methods start from this same initial point and use the same accurate eigensolver.
Standard HiSD and A-HiSD use a fixed step size of $2\times10^{-5}$, with momentum parameter $0.95$ for A-HiSD. HiSD-BB starts with a step size of $2\times10^{-5}$, with subsequent step sizes restricted to $[4\times10^{-6},10^{-4}]$.
The two limited-memory versions use $m=5$.
The fixed-step version takes $\beta_1=\beta_2=0.1$, while the full-memory and limited-memory line search versions use initial trial step sizes $\overline{\beta}_1=\overline{\beta}_2=1$ at each iteration.
The gradient norm tolerance is $10^{-8}$, and the maximum number of iterations is $500$.

\begin{figure}[!htbp]
    \centering
    \includegraphics[width=1\textwidth]{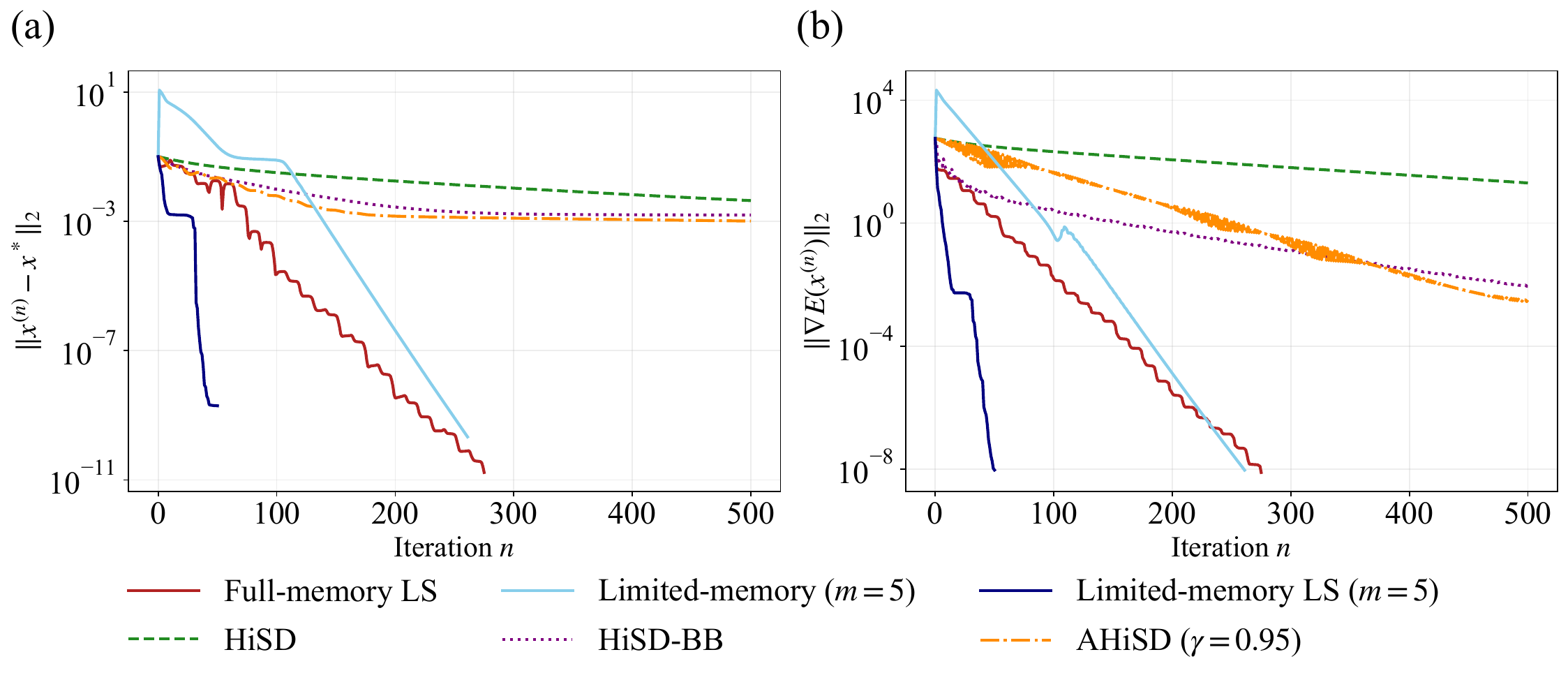}
    \caption{Plots of (a) error $\|x^{(n)}-x^*\|_2$ and
(b) gradient norm $\|\nabla E(x^{(n)})\|_2$ versus iteration number for different
saddle-search methods on the modified Rosenbrock example.}
    \label{fig:rosenbrock}
\end{figure}

Figure~\ref{fig:rosenbrock} shows the error
$\|x^{(n)}-x^*\|_2$ and gradient norm $\|\nabla E(x^{(n)})\|_2$ versus iteration number.
All three hybrid versions reach the prescribed tolerance within the iteration limit.
The limited-memory line search version performs best, requiring only about $50$ iterations.
The limited-memory fixed-step version exhibits a substantial initial increase in the error, whereas line search suppresses this increase and reaches the target accuracy sooner.
HiSD-BB and A-HiSD improve the performance of standard HiSD, but none of the three baseline methods reaches the prescribed tolerance within the iteration limit.
These results demonstrate the reduction in iteration counts achieved by the hybrid method on this ill-conditioned problem and the improvement provided by line search for the limited-memory version.

\subsection{Neural-network loss landscape: degenerate saddle manifolds}
We next consider a neural-network loss landscape to examine the performance of the proposed method near a degenerate saddle manifold. This experiment explores the method's numerical behavior beyond the nondegenerate setting of the local convergence analysis.

The fully connected neural network $f_\theta:\mathbb R^2\to\mathbb R$ in this experiment is 
\[
h_0(x)=x,\
h_\ell(x)=\tanh\!\left(W_\ell h_{\ell-1}(x)+b_\ell\right),
\ \ell=1,\ldots,4,
\quad
f_\theta(x)=w_5^\top h_4(x),
\]
where $\tanh$ acts componentwise,
$W_\ell\in\mathbb R^{d_\ell\times d_{\ell-1}}$,
$b_\ell\in\mathbb R^{d_\ell}$, and $w_5\in\mathbb R^{d_4}$.
We use layer widths $(d_0,\ldots,d_5)=(2,10,2,7,3,1)$, giving four hidden layers, a bias-free linear output layer, and $100$ trainable parameters.
The training dataset $\{x_i\}$ contains $N=400$ inputs sampled independently
and uniformly from $[-2,2]^2$, with target values
\[
y_i=f_{\mathrm {tar}}(x_i)+\frac{0.15}{c}r(x_i),
\quad
r(x)=\sum_{j=1}^{8}a_j\sin(\omega_j^\top x+\phi_j).
\]
Here, the fixed target network $f_{\mathrm {tar}}$ has the same form as $f_\theta$, but with widths $(2,2,2,2,3,1)$.
The parameters of $f_{\mathrm {tar}}$ and the coefficients $a_j$, $\omega_j$, $\phi_j$, and $c$ are fixed throughout the experiment.
The empirical loss serves as the energy function:
\[
E(\theta)=\frac1N\sum_{i=1}^{N}
\bigl(f_\theta(x_i)-y_i\bigr)^2.
\]

Following the embedding principle~\cite{zhang2021embedding}, we construct the reference saddle manifold $\mathcal M$ by embedding a local minimum of a narrower network $f_{\mathrm{nar}}$ into the wider network $f_\theta$ through neuron splitting as in~\cite{2026degenerate}.
The narrow network $f_{\mathrm {nar}}$ has layer widths $(2,2,2,2,3,1)$ and $30$ trainable parameters. It is trained on the same dataset $\{(x_i,y_i)\}_{i=1}^{N}$ by minimizing the empirical squared loss obtained from $E(\theta)$, with $f_\theta$ replaced by $f_{\mathrm{nar}}$.
The resulting saddle manifold $\mathcal M$ is a $13$-dimensional affine manifold.
At the reference point $\theta_*\in\mathcal M$, the Hessian has $31$ negative, $13$ zero, and $56$ positive eigenvalues.
The loss is approximately $3.71\times10^{-3}$, and the ratio between the largest and smallest nonzero eigenvalue magnitudes is approximately $3.46\times10^5$.

We compare standard HiSD, HiSD-BB, A-HiSD with momentum parameters $\gamma=0.9$ and $0.95$, and the full-memory and limited-memory line search versions of the proposed method.
The initial point is obtained by applying a random perturbation of Euclidean norm $0.01$ to $\theta_*$, and the same initial point is used for all methods.
Each method uses the $31$ accurate Hessian eigenpairs at every iteration.
Standard HiSD and both A-HiSD variants use a fixed step size of $0.1$. HiSD-BB starts with a step size of $0.1$, with subsequent step sizes restricted to $[0.02,0.5]$.
Both hybrid versions use the line search strategy described above, with initial trial step sizes $\overline\beta_1=\overline\beta_2=0.1$. The limited-memory version uses $m=10$. All methods are compared over $2000$ iterations.

\begin{figure}[!htbp]
    \centering
    \includegraphics[width=1\linewidth]{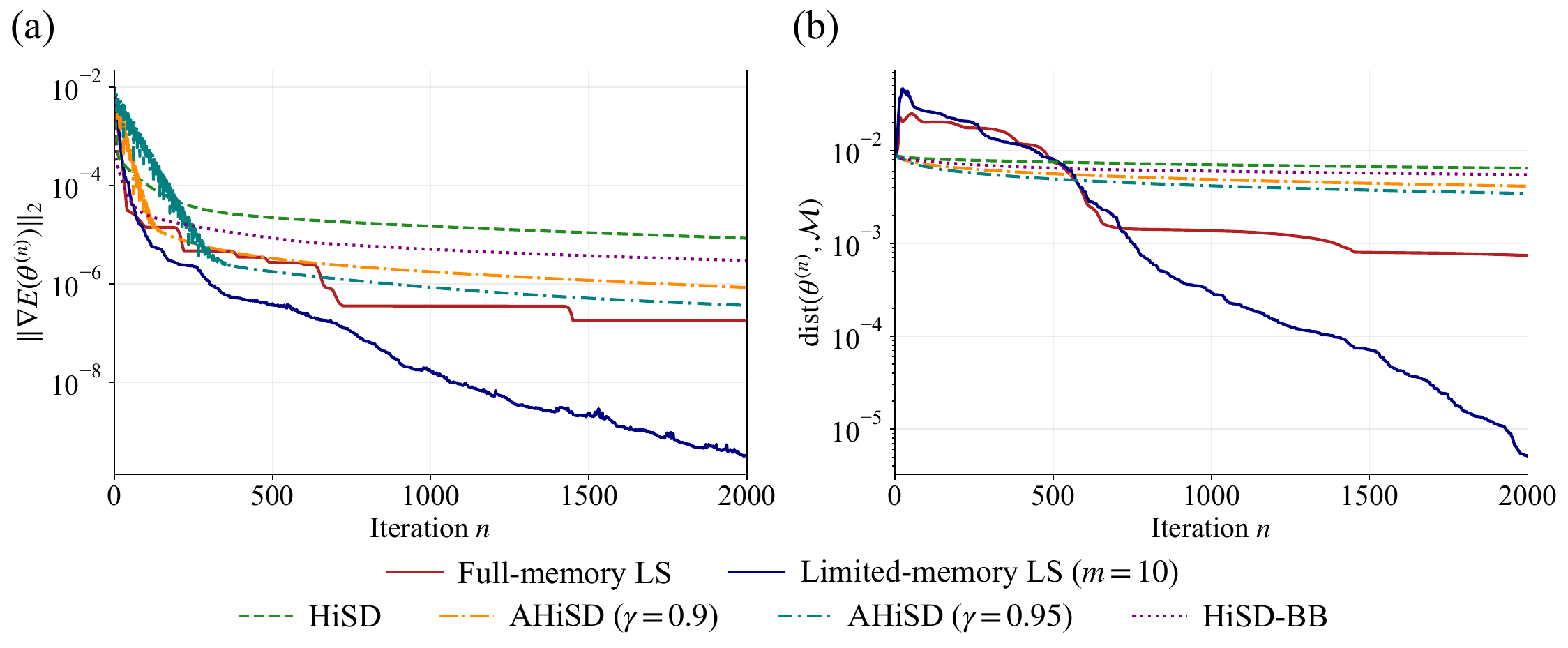}
    \caption{Plots of (a) gradient norm
    $\|\nabla E(\theta^{(n)})\|_2$ and
    (b) distance to the saddle manifold
    $\operatorname{dist}(\theta^{(n)},\mathcal M)$
    versus iteration number for different saddle-search methods.}
    \label{fig:neural-network}
\end{figure}

Figure~\ref{fig:neural-network} shows the gradient norm and the distance to the saddle manifold versus iteration number. Since $\mathcal M$ is affine, its orthogonal projection can be computed explicitly, allowing direct evaluation of this distance.
As illustrated in Figure~\ref{fig:neural-network}, BB step sizes and momentum acceleration improve the performance of standard HiSD. Both hybrid versions achieve higher accuracy within the same iteration budget.
The limited-memory line search version performs best, reducing the gradient norm and distance to the manifold to approximately $3\times10^{-10}$ and $5\times10^{-6}$, respectively. The full-memory line search version also attains smaller values of both quantities than the baseline methods. These results demonstrate the effectiveness of the proposed method in accelerating the local search for degenerate saddle points on critical manifolds.

\subsection{Landau--de Gennes model: upward and downward search}

Finally, we consider a Landau--de Gennes model with the bulk potential and material parameters used in~\cite{shi2023nonlinearity}. The energy is given by
\[
E[Q]
=
\int_{\Omega}
\left[
\frac12|\nabla Q|^2
+
\lambda
\left(
\frac{A}{4C}\operatorname{tr}(Q^2)
-\frac{B}{6C}\operatorname{tr}(Q^3)
+\frac18[\operatorname{tr}(Q^2)]^2
% -E_0
\right)
\right]\,dr,
\]
where $\Omega=[-1,1]^2$, $Q$ is a symmetric traceless $3\times 3$ tensor field
subject to prescribed Dirichlet boundary conditions.
In our experiment, we take
\[
\lambda=50,\quad
B=0.64\times10^4,\quad
C=0.35\times10^4,\quad
A=-\frac{B^2}{3C}.
\]
The domain is discretized on a uniform mesh with spacing $h=1/50$.
Discretizing the five independent tensor components gives $49005$ degrees of freedom.

We consider searches among the critical states $C+$, $M+$, and $N+$, using the state labels in~\cite{YinPRL}.
Their numerically computed Morse indices in the present model are $2$, $3$, and $4$, respectively.
Following the solution landscape framework in~\cite{yin2021searching}, we focus on upward and downward searches for the ill-conditioned saddle points $M+$ and $N+$.
The upward searches are $C+\to M+$ and $M+\to N+$, while the downward search is $N+\to M+$.
Each search starts from a known critical state with a perturbation of magnitude $10^{-3}$ along an appropriate Hessian eigenvector.
The prescribed index of HiSD is set to the index of the target saddle.

For this large-scale problem, we compare the limited-memory version of the proposed method with HiSD-BB and A-HiSD. HiSD-BB starts with a step size of $2\times10^{-4}$, with subsequent step sizes restricted to $[4\times10^{-5},10^{-3}]$.
A-HiSD uses a fixed step size of $8\times10^{-5}$ and a momentum parameter of $0.99$.
The limited-memory version uses $m=30$ and fixed step sizes
$\beta_1=\beta_2=10^{-2}$.
The eigenpairs are updated by a single step of the LOBPCG method in each iteration, and the iterations terminate when the gradient norm in the discrete $L^2$ metric is at most  $10^{-6}$.
Table~\ref{tab:ldg-iterations} reports the iteration counts and the wall time for the three searches.

\begin{table}[H]
\centering
\setlength{\tabcolsep}{3pt}
\caption{Number of iterations and wall time for various saddle searches at $\lambda=50$.}
\label{tab:ldg-iterations}
\begin{tabular}{cccccc}
\toprule
\multicolumn{2}{c}{Saddle search}
& \multirow{2}{*}{$k$}
& \multicolumn{3}{c}{Number of iterations (Wall time/s)} \\
\cmidrule(lr){1-2}
\cmidrule(lr){4-6}
From & To & & HiSD-BB & A-HiSD & Hybrid method \\
\midrule
$M+$ & $N+$ & 4
& 6,733,891 ($8.74\times10^{5}$)
& 287,484 (13,613)
& 2,604 (220) \\
$N+$ & $M+$ & 3
& 6,680,218 ($6.43\times10^{5}$)
& 287,320 (21,764)
& 2,473 (150) \\
$C+$ & $M+$ & 3
& 4,719,353 ($5.05\times10^{5}$)
& 206,425 (9,363)
& 4,430 (279) \\
\bottomrule
\end{tabular}
\end{table}

As shown in Table~\ref{tab:ldg-iterations},
the limited-memory hybrid method requires fewer iterations
in all three searches.
The improvement is particularly pronounced in the bidirectional
searches between $M+$ and $N+$, where the iteration counts
are reduced by approximately two orders of magnitude
compared with A-HiSD.
A substantial reduction is also observed for $C+\to M+$.
These results demonstrate the effectiveness of the limited-memory hybrid method in upward and downward saddle searches for this liquid-crystal model, with the reduction in iteration counts depending on the particular search.

\section{Concluding Remarks}
\label{sec:conclusions}

In this paper, we have proposed a hybrid Newton--BFGS HiSD method for accelerating saddle-point searches.
The proposed method preserves the ascent–descent structure of HiSD by combining Newton corrections on the unstable subspace with BFGS updates on the stable subspace. 
We apply projections to implement BFGS updates while the stable subspace is evolving during iterations.
Under suitable assumptions, our analysis establishes $Q$-superlinear convergence of the full-memory version.
We have also developed a limited-memory version and provided constant step size and line search strategies. 

Numerical experiments cover nonlinear test problems, neural-network loss landscapes, and the Landau–de Gennes model. Across these experiments, the proposed method demonstrates advantages in convergence speed and computational efficiency for ill-conditioned, degenerate, and large-scale saddle-point problems.
These results show that hybrid Newton--BFGS HiSD is a flexible method for both local refinement and upward and downward saddle searches, with favorable numerical performance and broad application potential.

An important direction for future work is to combine rapid convergence in each saddle point search with broader landscape exploration. Through this combination, we aim to develop methods that balance search efficiency with the breadth of exploration, allowing us to discover a wider variety of saddle points across different structural types and identify their connections. With richer solution landscapes, we can deepen our understanding of multistable structures and transition mechanisms in complex systems.
\appendix
% Format only this heading; keep reference and equation counters unchanged.
\makeatletter
\let\HiSDsavedseccntformat\@seccntformat
\def\@seccntformat#1{\appendixname~\csname the#1\endcsname.\hskip\betweenumberspace}
\section{Proofs of the Technical Lemmas}
\label{app:technical-proofs}
\let\@seccntformat\HiSDsavedseccntformat
\makeatother

Throughout this appendix, we use the standing setting, assumptions, notation,
and zero-extension convention of Section~\ref{sec:local-superlinear-convergence}.

\subsection{Proof of Lemma~\ref{lem:local-moving-estimates}}
\label{app:proof-local-moving-estimates}

Choose $R_0\le\min\{R,\mu/(4M)\}$. Assumption~\ref{ass:local-regularity} gives $\|H_{n+1}-H_n\|\le M\|s_n\|\le2MR_0\le\mu/2$. By the spectral separation and the subspace perturbation estimate in \cite[Theorem~V.3.6]{stewart1990matrix},
\begin{equation}\label{eq:app-projector-bound}
\|P_{n+1}^\pm-P_n^\pm\|_F
\le\sqrt d\,\|P_{n+1}^+-P_n^+\|
\le\frac{\sqrt d}{\mu-\mu/2}\|H_{n+1}-H_n\|
\le C_P\|s_n\|,
\end{equation}
where $C_P=2\sqrt d\,M/\mu$ and $P_n^-=I-P_n^+$.

To estimate the square roots, let $Z=S_{n+1}-S_n$. Since $P_{n+1}^-S_{n+1}=S_nP_n^-=0$, we have $P_{n+1}^-ZP_n^-=0$ and hence
\[
\|Z\|_F^2
=\|P_{n+1}^+Z\|_F^2+\|P_{n+1}^-Z\|_F^2
=\|P_{n+1}^+Z\|_F^2+\|P_{n+1}^-ZP_n^+\|_F^2
\le\|P_{n+1}^+Z\|_F^2+\|ZP_n^+\|_F^2.
\]
Since $S_{n+1}\succeq\sqrt\mu\,P_{n+1}^+$, we have
\[
\langle Z,S_{n+1}Z\rangle_F
=\operatorname{tr}(Z^\top S_{n+1}Z)
\ge\sqrt\mu\,\operatorname{tr}(Z^\top P_{n+1}^+Z)
=\sqrt\mu\,\|P_{n+1}^+Z\|_F^2.
\]
Similarly, $S_n\succeq\sqrt\mu\,P_n^+$ gives
$\langle Z,ZS_n\rangle_F\ge\sqrt\mu\,\|ZP_n^+\|_F^2$.

Combining these bounds with the preceding estimate and the identity
\[
S_{n+1}Z+ZS_n=S_{n+1}^2-S_n^2,
\]
we obtain
\[
\begin{aligned}
\sqrt\mu\,\|Z\|_F^2
&\le\sqrt\mu\bigl(\|P_{n+1}^+Z\|_F^2+\|ZP_n^+\|_F^2\bigr)\\
&\le\langle Z,S_{n+1}Z+ZS_n\rangle_F
\le\|Z\|_F\|S_{n+1}^2-S_n^2\|_F.
\end{aligned}
\]
Moreover, $S_j^2=H_jP_j^+$ gives
$S_{n+1}^2-S_n^2=(H_{n+1}-H_n)P_{n+1}^++H_n(P_{n+1}^+-P_n^+)$.
Together with \eqref{eq:app-projector-bound}, this yields
\begin{equation}\label{eq:app-square-root-bound}
\|S_{n+1}-S_n\|_F
\le\frac{\|S_{n+1}^2-S_n^2\|_F}{\sqrt\mu}
\le C_S\|s_n\|,
\end{equation}
where $C_S=(\sqrt d\,M+LC_P)/\sqrt\mu$. The first inequality also holds when $Z=0$. This proves \eqref{eq:projector-square-root-estimate} for $C_{\rm loc}\geq C_P+C_S$.

By the fundamental theorem of calculus and $P_{n+1}^+H_{n+1}=H_{n+1}P_{n+1}^+$,
\[
\eta_n-H_{n+1}\sigma_n
=P_{n+1}^+\int_0^1
\bigl(H(x^{(n)}+ts_n)-H_{n+1}\bigr)s_n\,dt.
\]
The segment joining the two iterates lies in $B_{R_0}(x^*)$. The norm of the right-hand side is therefore bounded by
$\int_0^1 M(1-t)\|s_n\|^2\,dt=\frac{M}{2}\|s_n\|^2$, proving \eqref{eq:projected-secant-estimate}.

Since $S_jP_j^+=P_j^+S_j=S_j$, direct expansion gives
\[
\widehat D_n-D_n
=(S_{n+1}-S_n)G_nS_{n+1}
+S_nG_n(S_{n+1}-S_n)
-(P_{n+1}^+-P_n^+).
\]
Using $\|S_j\|\le\sqrt L$ and \eqref{eq:app-projector-bound}--\eqref{eq:app-square-root-bound}, we obtain
\begin{equation}\label{eq:app-moving-space-bound}
\|\widehat D_n-D_n\|_F
\le\bigl(2\sqrt L\,C_S\|G_n\|+C_P\bigr)\|s_n\|
\le C_{\rm loc}(1+\|G_n\|)\|s_n\|.
\end{equation}
The choice $C_{\rm loc}=C_P+(1+2\sqrt L)C_S$ covers \eqref{eq:projector-square-root-estimate}--\eqref{eq:moving-potential-comparison}.

If Assumption~\ref{ass:projected-step} holds and $s_n\ne0$, then $\|\sigma_n\|\ge\chi\|s_n\|>0$. The secant estimate and the positive spectral bound imply
\[
\begin{aligned}
\sigma_n^\top\eta_n
&=\sigma_n^\top H_{n+1}\sigma_n
+\sigma_n^\top(\eta_n-H_{n+1}\sigma_n)\\
&\ge\mu\|\sigma_n\|^2-\tfrac{M}{2}\|\sigma_n\|\|s_n\|^2
\ge\left(\mu-\tfrac{M}{2\chi}\|s_n\|\right)\|\sigma_n\|^2.
\end{aligned}
\]
This proves \eqref{eq:local-curvature-lower-bound}, and the stated positive lower bound follows when $\|s_n\|\le\mu\chi/M$.

Finally, under the zero-extension convention,
\[
S(x^*)G_0S(x^*)-P_*^+
=S(x^*)\bigl[P_*^+G_0P_*^+-(H_*|_{V_*^\perp})^{-1}\bigr]S(x^*).
\]
Its Frobenius norm is at most $L\Delta_0^+$. The derivation of \eqref{eq:app-moving-space-bound} uses only the local bounds on $S$ and $P^+$, with the matrix $G_n$ held fixed, and therefore applies to any two points in $B_{R_0}(x^*)$. Applying this estimate to $x^{(0)}$ and $x^*$, with the same matrix $G_0$, and using the triangle inequality gives
\[
\begin{aligned}
\Phi_0
&\le\|S(x^*)G_0S(x^*)-P_*^+\|_F
+C_{\rm loc}(1+\|G_0\|)\|x^{(0)}-x^*\|\\
&\le L\Delta_0^++C_{\rm loc}(1+\Gamma_0)\|x^{(0)}-x^*\|,
\end{aligned}
\]
\nopagebreak
which proves \eqref{eq:initial-potential-comparison}.\qed

\subsection{Proof of Lemma~\ref{lem:bounded-deterioration}}
\label{app:proof-bounded-deterioration}

\emph{(i) Uniform boundedness of $G_n$.}
Every well-defined update preserves positive definiteness. By the definition of $D_j$ in \eqref{eq:current-positive-error},
\[
P_j^+G_jP_j^+=S_j^{-1}(P_j^++D_j)S_j^{-1},
\quad
\|P_j^+G_jP_j^+\|\le\frac{1+\overline\Phi}{\mu}.
\]
Set $b=(1+\overline\Phi)/\mu$ and $a_j=\|P_j^-G_jP_j^-\|$.
For any $z\in V_{j+1}$ with $\|z\|=1$, positive definiteness and the orthogonal decomposition at $x^{(j)}$ yield
\begin{equation}
\label{eq:quadratic-form}
\begin{aligned}
(z^\top G_jz)^{1/2}
&=\|G_j^{1/2}z\|
\le\sqrt{a_j}\,\|P_j^-z\|+\sqrt b\,\|P_j^+z\|\\
& = \sqrt{a_j}\,\|P_j^-z\|+\sqrt b\,\|(P_j^+-P_{j+1}^+)z\|
\le\sqrt{a_j}+\sqrt b\,\|P_{j+1}^--P_j^-\|,
\end{aligned}
\end{equation}
% where the last inequality uses $P_j^+z=(P_{j+1}^--P_j^-)z$.
Since $\sigma_j,\eta_j\in V_{j+1}^\perp$, multiplying the inverse-BFGS update on both sides by $P_{j+1}^-$ gives
\[
P_{j+1}^-G_{j+1}P_{j+1}^-=P_{j+1}^-G_jP_{j+1}^-.
\]
Taking the supremum over unit vectors $z\in V_{j+1}$ in the preceding quadratic-form estimate \eqref{eq:quadratic-form}, and using this identity and Lemma~\ref{lem:local-moving-estimates}, we obtain
\[
\sqrt{a_{j+1}}\le\sqrt{a_j}+C_{\rm loc}\sqrt b\,\|s_j\|.
\]
Since $a_0\le\Gamma_0$,
$\sqrt{a_j}
\le\sqrt{\Gamma_0}+C_{\rm loc}\sqrt b\sum_{\ell=0}^{j-1}\|s_\ell\|
\le\sqrt{\Gamma_0}+C_{\rm loc}\sqrt b\,\Lambda$.

For an arbitrary $z\in\mathbb R^d$, the same decomposition and the Cauchy--Schwarz inequality give
\[
\begin{aligned}
z^\top G_jz
&\le\bigl(\sqrt{a_j}\,\|P_j^-z\|+\sqrt b\,\|P_j^+z\|\bigr)^2\\
&\le(a_j+b)\bigl(\|P_j^-z\|^2+\|P_j^+z\|^2\bigr)
=(a_j+b)\|z\|^2.
\end{aligned}
\]
Thus $\|G_j\|\le a_j+b$, and \eqref{eq:finite-prefix-G-bound} follows with
$C_G=(\sqrt{\Gamma_0}+C_{\rm loc}\sqrt b\,\Lambda)^2+b$.

\emph{(ii) One-step bounded deterioration.}
Take $h=\mu\chi/M$. For $0<\|s_n\|\le h$, Lemma~\ref{lem:local-moving-estimates} gives
$\sigma_n^\top\eta_n\ge\frac{\mu}{2}\|\sigma_n\|^2>0$. Hence the inverse-BFGS update \eqref{eq:inverse-bfgs-update} is well defined and preserves symmetry and positive definiteness.

In this part only, write $P=P_{n+1}^+$ and $S=S_{n+1}$, and set $p=S\sigma_n$ and $q=S^{-1}\eta_n$. Then $p,q\in\operatorname{Ran}P$, $p^\top q=\sigma_n^\top\eta_n$, and $p/\|p\|=\omega_n$, as defined in \eqref{eq:weighted-secant-direction}. By Lemma~\ref{lem:local-moving-estimates} and the projected-step condition \eqref{eq:Projected-step-condition},
\[
\|q-p\|
=\|S^{-1}(\eta_n-H_{n+1}\sigma_n)\|
\le\frac{M}{2\sqrt\mu}\|s_n\|^2,\ 
\|p\|\ge\sqrt\mu\,\|\sigma_n\|\ge\sqrt\mu\,\chi\|s_n\|.
\]
Consequently, we have 
\begin{equation}
\label{eq:q-p-estimate}
\frac{\|q-p\|}{\|p\|}\le\frac{M}{2\mu\chi}\|s_n\|\le\frac12,
\quad
p^\top q\ge\|p\|^2-\|p\|\|q-p\|\ge\frac12\|p\|^2.
\end{equation}

The inverse-BFGS update \eqref{eq:inverse-bfgs-update} reads
\[
G_{n+1}
=\left(I-\frac{\sigma_n\eta_n^\top}{\sigma_n^\top\eta_n}\right)
G_n\left(I-\frac{\eta_n\sigma_n^\top}{\sigma_n^\top\eta_n}\right)
+\frac{\sigma_n\sigma_n^\top}{\sigma_n^\top\eta_n}.
\]
Using $S\sigma_n=p$, $Sq=\eta_n$, and $SP=PS=S$, we have
\[
\begin{aligned}
S\left(I-\frac{\sigma_n\eta_n^\top}{\sigma_n^\top\eta_n}\right)
&=S-\frac{p\eta_n^\top}{p^\top q}
=\left(P-\frac{pq^\top}{p^\top q}\right)S,\\
\left(I-\frac{\eta_n\sigma_n^\top}{\sigma_n^\top\eta_n}\right)S
&=S-\frac{\eta_np^\top}{p^\top q}
=S\left(P-\frac{qp^\top}{p^\top q}\right).
\end{aligned}
\]
Multiplying the update \eqref{eq:inverse-bfgs-update} on both sides by $S$ and using $SG_nS=P+\widehat D_n$ and $SG_{n+1}S=P+D_{n+1}$ gives
\begin{equation}\label{eq:app-weighted-bfgs}
P+D_{n+1}
=\left(P-\frac{pq^\top}{p^\top q}\right)
(P+\widehat D_n)
\left(P-\frac{qp^\top}{p^\top q}\right)
+\frac{pp^\top}{p^\top q}.
\end{equation}

The projected secant estimate \eqref{eq:projected-secant-estimate} shows that $\eta_n$ and $H_{n+1}\sigma_n$ are close, so we first consider the simpler case of equality and then estimate the effect of their difference.
If $\eta_n=H_{n+1}\sigma_n$ in \eqref{eq:inverse-bfgs-update}, then we have $q=S^{-1}H_{n+1}\sigma_n=p$
and the right hand side of \eqref{eq:app-weighted-bfgs} becomes
$(P-\omega_n\omega_n^\top)(P+\widehat D_n)
(P-\omega_n\omega_n^\top)+\omega_n\omega_n^\top$.
Generally, we have the identity
\begin{equation}\label{eq:app-ideal-weighted-error}
\begin{aligned}
(P-\omega_n\omega_n^\top)(P+\widehat D_n)
(P-\omega_n\omega_n^\top)+\omega_n\omega_n^\top
=P+
(P-\omega_n\omega_n^\top)\widehat D_n(P-\omega_n\omega_n^\top),
\end{aligned}
\end{equation}
where we use that $P-\omega_n\omega_n^\top$ is an orthogonal projector. 
Employing the symmetry of $\widehat D_n$ and
$P\widehat D_n=\widehat D_nP=\widehat D_n$, we can estimate
$\|(P-\omega_n\omega_n^\top)\widehat D_n(P-\omega_n\omega_n^\top)\|_F$ as follows:
\[
\begin{aligned}
&\|(P-\omega_n\omega_n^\top)\widehat D_n(P-\omega_n\omega_n^\top)\|_F^2
\le \|\widehat D_n(P-\omega_n\omega_n^\top)\|_F^2\\
=\;& \operatorname{tr}\left((P-\omega_n\omega_n^\top){\widehat D_n^2} (P-\omega_n\omega_n^\top)\right)=\operatorname{tr}\left({\widehat D_n^2} (P-\omega_n\omega_n^\top)\right)\\
=\;&\operatorname{tr}(\widehat D_n^2)-\omega_n^\top\widehat D_n^2\omega_n
=\widehat\Phi_n^2-\|\widehat D_n\omega_n\|^2.
\end{aligned}
\]
Therefore, we have $\|(P-\omega_n\omega_n^\top)\widehat D_n(P-\omega_n\omega_n^\top)\|_F\le \widehat\Phi_n\le\overline\Phi<1$, and
\[
\begin{aligned}
\|\widehat D_n\omega_n\|^2
&\le\widehat\Phi_n^2-
\bigl\|(P-\omega_n\omega_n^\top)\widehat D_n(P-\omega_n\omega_n^\top)\bigr\|_F^2\\
&\le2\Bigl(\widehat\Phi_n-
\bigl\|(P-\omega_n\omega_n^\top)\widehat D_n(P-\omega_n\omega_n^\top)\bigr\|_F\Bigr).
\end{aligned}
\]
Finally, we have
\begin{equation}\label{eq:app-ideal-descent}
\bigl\|(P-\omega_n\omega_n^\top)\widehat D_n(P-\omega_n\omega_n^\top)\bigr\|_F
\le\widehat\Phi_n-\tfrac12\|\widehat D_n\omega_n\|^2.
\end{equation}

Subtracting \eqref{eq:app-ideal-weighted-error} from
\eqref{eq:app-weighted-bfgs}, we obtain the following decomposition:
\[
\begin{aligned}
&\ D_{n+1}-(P-\omega_n\omega_n^\top)\widehat D_n(P-\omega_n\omega_n^\top)
=\underbrace{\left(\omega_n\omega_n^\top-\frac{pq^\top}{p^\top q}\right)
(P+\widehat D_n)\left(P-\frac{qp^\top}{p^\top q}\right)}_{R_1}\\
&\qquad\qquad\qquad\qquad\qquad
+\underbrace{(P-\omega_n\omega_n^\top)(P+\widehat D_n)
\left(\omega_n\omega_n^\top-\frac{qp^\top}{p^\top q}\right)}_{R_2}
+\underbrace{\left(\frac{pp^\top}{p^\top q}-\omega_n\omega_n^\top\right)}_{R_3}.
\end{aligned}
\]
We now estimate the Frobenius norms of $R_1$, $R_2$, and $R_3$ separately.
Since $\omega_n=p/\|p\|$ and $p,q\in\operatorname{Ran}P$,
we have $(P-\omega_n\omega_n^\top)p=0$ and
\[
\frac{qp^\top}{p^\top q}-\omega_n\omega_n^\top
=\frac{(P-\omega_n\omega_n^\top)qp^\top}{p^\top q}
=\frac{(P-\omega_n\omega_n^\top)(q-p)p^\top}{p^\top q}.
\]
Using \(p^\top q\ge \frac12\|p\|^2\), we obtain
\[
\left\|\frac{qp^\top}{p^\top q}-\omega_n\omega_n^\top\right\|
\le\left\|\frac{qp^\top}{p^\top q}-\omega_n\omega_n^\top\right\|_F
\le \frac{2\|q-p\|}{\|p\|}.
\]
Since \(P-\omega_n\omega_n^\top\) is an orthogonal projector, the triangle inequality then gives
\[
\left\|P-\frac{qp^\top}{p^\top q}\right\|\le \|P-\omega_n\omega_n^\top\|+\left\|\frac{qp^\top}{p^\top q}-\omega_n\omega_n^\top\right\|
\le 1+\frac{2\|q-p\|}{\|p\|}
\le 2.
\]
Notice that $\|P+\widehat D_n\|\le1+\overline\Phi\le2$, so we have
$$
\|R_1\|_F\le4
\left\|\omega_n\omega_n^\top-\frac{pq^\top}{p^\top q}\right\|_F\le\frac{8\|q-p\|}{\|p\|};\ \|R_2\|_F
\le2\left\|\omega_n\omega_n^\top-\frac{qp^\top}{p^\top q}\right\|_F
\le\frac{4\|q-p\|}{\|p\|}.
$$
Substituting $\omega_n=p/\|p\|$, we obtain
\[
\|R_3\|_F=\left\|\frac{p^\top(p-q)}{p^\top q}\frac{pp^\top}{p^\top p}\right\|_F
\le \frac{|p^\top(p-q)|}{p^\top q}
\le2\frac{\|q-p\|}{\|p\|}.
\]
Combining the above bounds and \eqref{eq:q-p-estimate}, we obtain the following estimate:
\begin{equation}\label{eq:app-bfgs-perturbation}
\bigl\|D_{n+1}-(P-\omega_n\omega_n^\top)\widehat D_n(P-\omega_n\omega_n^\top)\bigr\|_F
\le14\frac{\|q-p\|}{\|p\|}
\le\frac{7M}{\mu\chi}\|s_n\|.
\end{equation}

Applying the triangle inequality and using \eqref{eq:app-ideal-descent} and \eqref{eq:app-bfgs-perturbation}, we obtain
\[
\Phi_{n+1}
\le\widehat\Phi_n-\tfrac12\|\widehat D_n\omega_n\|^2+\frac{7M}{\mu\chi}\|s_n\|.
\]
Finally, Lemma~\ref{lem:local-moving-estimates} and $\|G_n\|\le\overline\Gamma$ give
\[
\widehat\Phi_n\le\Phi_n+\|\widehat D_n-D_n\|_F
\le\Phi_n+C_{\rm loc}(1+\overline\Gamma)\|s_n\|.
\]
This proves \eqref{eq:bfgs-deterioration-current-space} with $c_{\rm B}=1/2$ and
$C_{\rm B}=C_{\rm loc}(1+\overline\Gamma)+7M/(\mu\chi)$.\qed

\subsection{Proof of Lemma~\ref{lem:inverse-dennis-more}}
\label{app:proof-inverse-dennis-more}

Under the standing setting \eqref{ass:123}, every well-defined update preserves $G_n=G_n^\top\succ0$. If $s_n=0$, the split step \eqref{eq:local-split-step} gives
\[
0=P_n^-s_n=-(H_n|_{V_n})^{-1}P_n^-g_n, 
\qquad
0=P_n^+s_n=-P_n^+G_nP_n^+g_n.
\]
Since $H_n|_{V_n}$ is invertible and $P_n^+G_nP_n^+$ is positive definite on $V_n^\perp$, it follows that
\[
P_n^-g_n=P_n^+g_n=0,\quad g_n=0,
\]
contradicting nontermination, which indicates $s_n\ne0$. The projected-step condition \eqref{eq:Projected-step-condition} then gives $\|\sigma_n\|\ge\chi\|s_n\|>0$, so $\omega_n$ is well defined.

From $\sum_{n=0}^{\infty}\|s_n\|<\infty$, we have $\|s_n\|\to 0$, as $n\rightarrow \infty$. Let $h,c_{\rm B},C_{\rm B}$ be the constants in Lemma~\ref{lem:bounded-deterioration}(ii) corresponding to $\overline\Gamma,\overline\Phi$, and choose $n_0$ such that $0<\|s_n\|\le h$ for all $n\ge n_0$. By Lemma \ref{lem:bounded-deterioration}, for every $N\ge n_0$,
\[
c_{\rm B}\sum_{n=n_0}^N\|\widehat D_n\omega_n\|^2
\le\Phi_{n_0}-\Phi_{N+1}+C_{\rm B}\sum_{n=n_0}^N\|s_n\|
\le\Phi_{n_0}+C_{\rm B}\sum_{n=n_0}^{\infty}\|s_n\|<\infty.
\]
The right-hand side is independent of $N$, hence $\|\widehat D_n\omega_n\|\to0$, as $n\rightarrow \infty$. Recalling \eqref{eq:weighted-secant-direction},
$\omega_n=\frac{S_{n+1}\sigma_n}{\|S_{n+1}\sigma_n\|}$,
and using $S_{n+1}s_n=S_{n+1}\sigma_n$ and $\|S_{n+1}\sigma_n\|=\|S_{n+1}s_n\|\le\sqrt L\,\|s_n\|$, we obtain
\begin{equation}\label{eq:app-new-space-direction}
\frac{\|\widehat D_nS_{n+1}s_n\|}{\|s_n\|}
=\frac{\|\widehat D_n\omega_n\|\|S_{n+1}\sigma_n\|}{\|s_n\|}
\le\sqrt L\,\|\widehat D_n\omega_n\|\longrightarrow0, \quad\text{as} \  n\rightarrow \infty.
\end{equation}

The estimate in \eqref{eq:app-new-space-direction} concerns the positive subspace at $x^{(n+1)}$,
whereas \eqref{eq:current-space-inverse-dennis-more} concerns
the positive subspace at $x^{(n)}$. To relate the two, write
\[
D_nS_ns_n-\widehat D_nS_{n+1}s_n
=(D_n-\widehat D_n)S_ns_n+\widehat D_n(S_n-S_{n+1})s_n.
\]
By Lemma~\ref{lem:local-moving-estimates}, $\|G_n\|\le\overline\Gamma$, and $\|\widehat D_n\|\le\overline\Phi$, we have
\[
\begin{aligned}
\|D_nS_ns_n-\widehat D_nS_{n+1}s_n\|
&\le\bigl(\sqrt L\,\|D_n-\widehat D_n\|_F
+\overline\Phi\,\|S_n-S_{n+1}\|_F\bigr)\|s_n\|\\
&\le C_{\rm loc}\bigl((1+\overline\Gamma)\sqrt L+\overline\Phi\bigr)\|s_n\|^2.
\end{aligned}
\]
Combining the preceding estimate with \eqref{eq:app-new-space-direction} and using $\|s_n\|\to0$, as $n\rightarrow \infty$, we obtain
\begin{equation}\label{eq:app-current-weighted-direction}
\frac{\|D_nS_ns_n\|}{\|s_n\|}\longrightarrow0, \quad\text{as} \  n\rightarrow \infty.
\end{equation}
By the definition of $D_n$ and the zero-extension convention,
\[
S_n^{-1}D_nS_n
=P_n^+G_nP_n^+S_n^2-P_n^+
=P_n^+G_nP_n^+H_nP_n^+-P_n^+.
\]
Using $S_nP_n^+=S_n$ and $\|S_n^{-1}|_{V_n^{\perp}}\|\le\mu^{-1/2}$, it follows from \eqref{eq:app-current-weighted-direction} that
\[
\begin{aligned}
\quad\frac{\bigl\|(P_n^+G_nP_n^+H_nP_n^+-P_n^+)P_n^+s_n\bigr\|}{\|s_n\|}
&=\frac{\|S_n^{-1}D_nS_ns_n\|}{\|s_n\|}\\
&\le\frac1{\sqrt\mu}\frac{\|D_nS_ns_n\|}{\|s_n\|}\longrightarrow0,\quad\text{as} \  n\rightarrow \infty,
\end{aligned}
\]
which proves \eqref{eq:current-space-inverse-dennis-more}.

Finally, set $K_n=P_n^+G_nP_n^+$. Recalling that $\mathcal R_n=g_n+H_ns_n$ and using the split step \eqref{eq:local-split-step}, we obtain
\[
P_n^-\mathcal R_n
=P_n^-g_n+H_nP_n^-s_n
=P_n^-g_n-H_n(H_n|_{V_n})^{-1}P_n^-g_n=0,
\]
which indicates $\mathcal R_n\in V_n^\perp$. 
Moreover, $P_n^+s_n=-K_ng_n$ and $H_nP_n^+=P_n^+H_n$ imply
\[
K_n\mathcal R_n
=K_ng_n+K_nH_ns_n
=-P_n^+s_n+K_nH_nP_n^+s_n
=(K_nH_nP_n^+-P_n^+)P_n^+s_n.
\]
By \eqref{eq:current-space-inverse-dennis-more}, this identity gives $\|K_n\mathcal R_n\|/\|s_n\|\to0$.
Since $\mathcal R_n\in V_n^\perp$, it remains to establish a uniform positive lower bound for $K_n$ on $V_n^\perp$ to obtain the desired residual estimate \eqref{eq:newton-residual}.

Since $K_n=S_n^{-1}(P_n^++D_n)S_n^{-1}$, for every $z\in V_n^\perp$, the bound $\|D_n\|\le\overline\Phi<1$ and the upper spectral bound give
\[
z^\top K_nz
=(S_n^{-1}z)^\top(P_n^++D_n)(S_n^{-1}z)
\ge(1-\overline\Phi)\|S_n^{-1}z\|^2
\ge\frac{1-\overline\Phi}{L}\|z\|^2,
\]
which indicates $\|K_nz\|\ge(1-\overline\Phi)\|z\|/L$.
By taking $z=\mathcal R_n\in V_n^\perp$ and using \eqref{eq:current-space-inverse-dennis-more}, we obtain
\[
\begin{aligned}
\frac{\|\mathcal R_n\|}{\|s_n\|}
&\le\frac{L}{1-\overline\Phi}\frac{\|K_n\mathcal R_n\|}{\|s_n\|}\\
&=\frac{L}{1-\overline\Phi}
\frac{\bigl\|(K_nH_nP_n^+-P_n^+)P_n^+s_n\bigr\|}{\|s_n\|}
\longrightarrow0, \quad\text{as} \  n\rightarrow \infty.
\end{aligned}
\]
Together with $P_n^-\mathcal R_n=0$, this proves \eqref{eq:newton-residual}.\qed

\section*{Declarations}

\begingroup
\small
\noindent{\sffamily\bfseries Funding}\enspace
This work was supported by the National Natural Science Foundation of China (No. 12225102, T2321001, 12288101, 12601819).
\par\medskip

\noindent{\sffamily\bfseries Code and Data Availability}\enspace
The code and accompanying input data for the numerical experiments presented in this paper are publicly available at
\url{https://github.com/omega31415/hybrid-newton-bfgs-hisd}.
\par\medskip

\noindent{\sffamily\bfseries Conflict of interest}\enspace
The authors have no competing interests to declare that are
relevant to the content of this article.
\par

\endgroup

\bibliographystyle{spmpsci}
\bibliography{refs}

\end{document}